\documentclass[reqno,fleqn,11pt]{amsart}%
\usepackage{ifthen} 

\provideboolean{shownotes} 
\setboolean{shownotes}{true} 
\newcommand{\margnote}[1]{
\ifthenelse{\boolean{shownotes}}%
{\marginpar{\raggedright\tiny\texttt{#1}}}%
{}%
}
\newcommand{\hole}[1]{
\ifthenelse{\boolean{shownotes}}%
{\begin{center} \fbox{ \rule {.25cm}{0cm}
\rule[-.1cm]{0cm}{.4cm} \parbox{.85\textwidth}{\begin{center}
\texttt{#1}\end{center}} \rule {.25cm}{0cm}}\end{center}}
{}
}
\newtheorem{theorem}{Theorem}[section]
\newtheorem{proposition}[theorem]{Proposition}
\newtheorem{lemma}[theorem]{Lemma}

\theoremstyle{remark}
\newtheorem{remark}[theorem]{Remark}
\newtheorem{definition}[theorem]{Definition}

\newcommand{\N}{\mathbb{N}}
\newcommand{\R}{\mathbb{R}}
\newcommand{\e}{\varepsilon}
\newcommand{\bue}{\bar u^{\e}}
\newcommand{\bve}{\bar v^{\e}_{2}}

\numberwithin{equation}{section}

\begin{document}

\title [Shock waves for radiative hyperbolic--elliptic systems]
{Shock waves for radiative\\ hyperbolic--elliptic systems}

\author{Corrado Lattanzio}
\address{Corrado Lattanzio --- Sezione di Matematica per l'Ingegneria\\
Dipartimento di Matematica Pura ed Applicata \\
Universit\`a di L'Aquila \\
Piazzale E. Pontieri, 2 \\
		     Monteluco di Roio\\
		     67040 L'Aquila, Italy}
\email{corrado@univaq.it}

\author{Corrado Mascia}
\address{Corrado Mascia --- 
Dipartimento di Matematica ``G. Castelnuovo''\\
Universit\`a di Roma ``La Sapienza''\\
Piazzale A. Moro, 2 \\
		     00185 Roma, Italy}
\email{mascia@mat.uniroma1.it}

\author{Denis Serre}
\address{Denis Serre --- Unit\'e de Math\'ematiques Pures et Appliqu\'ees  \\ 
UMR CNRS 5669 \\ ENS Lyon \\
46, All\'ee d'Italie \\ 69364 Lyon Cedex 07, France}
\email{Denis.SERRE@umpa.ens-lyon.fr}

\baselineskip 14pt

\begin{abstract}
The present paper deals with the following hyperbolic--elliptic coupled system, 
modelling dynamics of a gas in presence of radiation, 
\begin{equation*}
    \begin{cases}
	u_{t}+ f(u)_{x} +Lq_{x}=0 & \\
	-q_{xx} + Rq +G\cdot u_{x} =0, &
    \end{cases}\qquad x\in\R,\quad t>0,
\end{equation*}
where $u\in\R^{n}$, $q\in\R$ and $R>0$, $G$, $L\in\R^{n}$. 
The flux function $f\,:\,\R^n\to\R^n$ is smooth and
such that $\nabla f$ has $n$ distinct real eigenvalues for any $u$.
 
The problem of existence of {\it admissible radiative shock wave} is considered,
i.e. existence of a solution of the form $(u,q)(x,t):=(U,Q)(x-st)$, such that 
$(U,Q)(\pm\infty)=(u_\pm,0)$, and $u_\pm\in\R^n$, $s\in\R$ define a shock wave for the 
reduced hyperbolic system, obtained by formally putting $L=0$. 

It is proved that, if $u_-$ is such that $\nabla \lambda_{k}(u_-) \cdot r_{k}(u_-)\neq 0$,
(where $\lambda_k$ denotes the $k$-th eigenvalue of $\nabla f$ and $r_k$ a 
corresponding right eigenvector) and
\begin{equation*}
   (\ell_{k}(u_{-})\cdot L)\,(G\cdot r_{k}(u_{-})) >0,
\end{equation*}
then there exists a neighborhood $\mathcal U$ of $u_-$ such that
for any $u_+\in{\mathcal U}$, $s\in\R$ such that the triple $(u_{-},u_{+};s)$ 
defines a shock wave for the reduced hyperbolic system, 
there exists a (unique up to shift) admissible radiative shock wave for
the complete hyperbolic--elliptic system.
The proof is based on reducing the system case to the scalar case, 
hence the problem of existence for the scalar case with general strictly convex fluxes
is considered, generalizing existing results for the Burgers' flux $f(u)=u^2/2$.
Additionally, we are able to prove that the profile $(U,Q)$  gains 
smoothness when the size of the shock $|u_+-u_-|$ is small enough, as
previously proved for the Burgers' flux case.

Finally, the general case of nonconvex fluxes is also treated, showing similar 
results of existence and regularity for the profiles.
\end{abstract}

\date{}
\maketitle

\thispagestyle{empty}

\section{Introduction}\label{sec:intro}
The dynamics of a gas in presence of radiation, due to high--temperature effects, can be 
modeled by compressible Euler equations with an additional term in the flux of energy.
Dealing with small perturbation of a fixed equilibrium state in one space dimension, 
this leads to consider an hyperbolic--elliptic coupled system of the form
\begin{equation}
    \begin{cases}
	u_{t}+ f(u)_{x} +Lq_{x}=0 & \\
	-q_{xx} + Rq +G\cdot u_{x} =0, &
    \end{cases}
    \label{eq:systemintro}
\end{equation}
where $x\in\R$, $t>0$, $u\in\R^{n}$, $q\in\R$ and $R>0$, 
$G$, $L\in\R^{n}$ are constant vectors. 
The flux function $f\,:\,\R^n\to\R^n$ is assumed to be smooth and
such that the reduced system
\begin{equation}
   u_{t}+ f(u)_{x}=0
\label{eq:reduced}
\end{equation}
is strictly hyperbolic, i.e. $\nabla f(u)$ has $n$ distinct real eigenvalues
for any state $u$ under con\-si\-de\-ra\-tion.
For later use, we denote such eigenvalues with $\lambda_1(u)<\dots<\lambda_n(u)$
and with $\ell_1(u),\dots, \ell_n(u)$, $r_1(u),\dots, r_n(u)$ the corresponding
left and right eigenvectors normalized so that $\ell_i(u)\,\cdot\,r_j(u)=\delta_{ij}$ for
any $i,j$.

The system (\ref{eq:systemintro}) can be obtained from the complete gas 
dynamics equation with heat--flux radiative term by using a 
differential approximation of the integral equation for the
radiative term.
This approach has been proposed in the pioneering paper \cite{Ham71}. 
After that, equations (\ref{eq:systemintro}) are also called {\it Hamer system
for radiating gas}
(see also \cite{KN01,KT04} for details on the derivation).
Since the quantity $q$ represents the radiative heat-flux term,
considering it as a scalar quantity is physically meaningful
(see \cite{VK65} for a complete physical description of the phenomenon).

Under the hyperbolic rescaling
$(\partial_t,\partial_x)\;\mapsto\;(\varepsilon\,\partial_t,\varepsilon\,\partial_x)$,
the system (\ref{eq:systemintro}) becomes
\begin{equation*}
    \begin{cases}
	u_{t}+ f(u)_{x} +Lq_{x}=0 & \\
	-\varepsilon^2\,q_{xx} + Rq +\varepsilon\,G\cdot u_{x} =0. &
    \end{cases}
\end{equation*}
Eliminating the $q$ variable, 
 \begin{equation}
   u_{t}+ f(u)_{x}
   =\varepsilon\,R^{-1}L\otimes G\,u_{xx}
     +\varepsilon^2\,R^{-1}\bigl(u_{t}+ f(u)_{x}\bigr)_{xx}.
 \label{eq:expansion}
 \end{equation}
Hence the rescaled hyperbolic--elliptic system can be rewritten as
a singular perturbation of a system of conservation laws. 
In particular, for $\varepsilon$ sufficiently small, the system can be seen 
as a higher order correction of a viscous system of conservation laws
with (degenerate) diffusion term given by the rank--one 
matrix $\varepsilon\,R^{-1}L\otimes G$.
In particular, this suggests that qualitative properties of solutions of
(\ref{eq:systemintro}) should resemble analogous properties of viscous 
system of conservation law as soon as:
$u$ varies along a direction not orthogonal to both vectors $L$ and $G$ {\it (non degeneracy)};
variations of $u$ are mainly in the small frequencies regime  {\it (small perturbations)}.
Indeed, failing of the first condition would imply degeneration of the diffusion term 
$\varepsilon\,R^{-1}L\otimes G\,u_{xx}$, and failing  of the second would give to the 
higher order term $\varepsilon^2\,R^{-1}\bigl(u_{t}+ f(u)_{x}\bigr)_{xx}$ a dominating r\^ole.

The above effects are present also in the {\it scalar case}, i.e. $u\in\R$
(for a complete and introductory presentation see \cite{Ser04}).
As noted in \cite{ST92}, in this situation, the $2\times 2$ system (\ref{eq:systemintro}) 
enjoys many of the properties of scalar viscous conservation laws:
$L^1$--contraction, comparison principle, conservation of mass,
constant solutions. 
Thanks to these properties, global existence and uniqueness of solutions
can be proved (see \cite{Ito97} for data in $BV$, \cite{LM03}  for
data in $L^1\cap L^\infty$, \cite{DiF06} for the multidimensional case 
for data in $L^1\cap L^\infty$).
Nevertheless, regularization property does not hold:
in \cite{KN99b} and in \cite{LT01} (with a more detailed description)
it is shown that there are initial data such that the corresponding solution to
the Cauchy problems develops discontinuity in finite time.

The loss of regularity appears also when dealing with the 
problem\footnote{The analogous problem for the viscous regularization of a system
of hyperbolic conservation laws is sometimes referred to as {\it Gel'fand problem},
and the first rigorous mathematical result has been proved in \cite{Gil51}.}

\begin{quote}\it
given $u_\pm$, asymptotic states of an admissible shock wave solution to (\ref{eq:reduced}),
does there exist a traveling wave solution to (\ref{eq:systemintro}) with same speed
of the shock and asymptotic states $(u_\pm,0)$?
\end{quote}

From now on, we refer to such a solution as a {\sc radiative shock wave}.
In the scalar case and for $f(s)=\frac12\,s^2$, in \cite{KN99} it has proved that the answer is 
affirmative, but that the profile of the radiative shock wave is discontinuous whenever the 
hyperbolic shock is large enough, i.e. $|u_--u_+|$ is large. 
The precise statement will be recalled later on.
In the same article, stability and decay rate of perturbations are determined.
The absence/presence of jumps in small/large radiative shock wave is again a 
manifestation of regularity properties for solution of (\ref{eq:systemintro}): 
small transitions can be obtained through smooth solution, large transitions cannot.

For the scalar case, many other results are available in literature.
For the sake of completeness, let us mention them, collected in two different groups:\\
-- {\it Large--time behavior}: to prove stability, possibly with decay rate of the perturbations, of 
constant states \cite{Ito97,Ser03}, shock profiles \cite{KN99,Ser04},
rarefaction waves \cite{KT04},
and to find asymptotic profiles for such perturbations \cite{Lau05,DiFL05};\\
-- {\it Weak solutions}:  to find evolution/regularity of discontinuity curves \cite{Nis00},
to determine relaxation limit under hyperbolic/parabolic rescaling \cite{DiF06, LM03}.
\medskip

The theory for system is still at the very beginning and just few results are available.
The first paper in this direction is \cite{KNN98}, where global
existence, asymptotic behavior and decay rate are proved for the solution
to the Cauchy problem with initial data that are small perturbations of constant
states.
Generalizations have been given in \cite{IK02} and \cite{KNN03}, especially
in the precise description of asymptotic profiles (diffusion waves).
In particular, it has been proved that, for large time, the solutions to (\ref{eq:systemintro})
are well--approximated by the solutions to the viscous system of conservation 
laws obtained from (\ref{eq:expansion}) when disregarding the $O(\varepsilon^2)$ term.
The singular limit counter part, i.e. $\varepsilon\to 0^+$, has been dealt with in
\cite{KN01}.
A different approach for analyzing the singular limit, based on the notion of positively 
invariant domain, has been considered in \cite{Ser04b}.
\medskip

The present paper deals with the problem of proving existence of radiative shock waves
in the case of general systems of the form (\ref{eq:systemintro}).
First of all let us recall the definition of shock wave and radiative shock wave,
where, for the sake of clarity, we refer here to the genuinely nonlinear case.

\begin{definition}\label{def:shock}
A {\sc shock wave} of the hyperbolic system (\ref{eq:reduced}) is
a weak solution of the form
 \begin{equation*}
  u(x,t):=u_-\chi_{{}_{(-\infty,x_0)}}(x-st)+u_+\chi_{{}_{(x_0,+\infty)}}(x-st),
  \qquad\qquad x_0\in\R
 \end{equation*} 
where $u_\pm\in\R^n$, $s\in\R$ satisfy, for some $k\in\{1,\dots,n\}$, the conditions
$\lambda_{k-1}(u_-)<s<\lambda_k(u_-)$ and 
$\lambda_k(u_+)<s<\lambda_{k+1}(u_+)$,
where $\lambda_1(u)<\dots<\lambda_n(u)$ denote the (real) eigenvalues of $\nabla f(u)$
and $\chi_{{}_{I}}(x)$ is the characteristic function of the set $I$.
\end{definition}

In the general case, the entropy condition for Lax shocks, $\lambda_{k-1}(u_-)<s<\lambda_k(u_-)$ and $\lambda_k(u_+)<s<\lambda_{k+1}(u_+)$, should be replaced by the Liu--E condition
(see the end of Section \ref{sec:system} for details).

\begin{definition}\label{def:radiative}
    A {\sc radiative shock wave} of the hyperbolic--elliptic system (\ref{eq:systemintro}) is
a weak solution $(u,q)(x,t):=(U,Q)(x-st)$ such that
 \begin{equation*}
   \lim_{\xi\to\pm\infty} (U,Q)(\xi)=(u_\pm,0),
 \end{equation*}
where $u_\pm\in\R^n$, $s\in\R$ defines a shock wave for the reduced hyperbolic system 
(\ref{eq:reduced}).
\end{definition}

The usual Rankine--Hugoniot condition 
\begin{equation}
  f(u_+)-f(u_-)=s(u_+-u_-),
 \label{eq:rh}
 \end{equation}
relating the states $u_\pm$ and the speed of propagation $s$, holds also for radiative shock waves.
This is readily seen by integrating over all $\R$ the first equation (\ref{eq:systemintro}) 
and taking in account the asymptotic limits of the wave.
 
Being a weak solution, a radiative shock wave may be discontinuous at some $\xi_0$.
Let $(U,Q)$ be piecewise $C^1$ in a neighborhood $(\xi_-,\xi_+)$ of $\xi_0$, that is
 \begin{equation*}
   U\in C^1((\xi_-,\xi_0])\cap C^1([\xi_0,\xi_+)),\qquad
   Q\in C^2((\xi_-,\xi_0])\cap C^2([\xi_0,\xi_+)).
 \end{equation*}
If the profile $U(\xi_0-)\neq U(\xi_0+)$, then $U'$ has a delta term concentrated at
$\xi_0$ and, as a consequence of the second equation in (\ref{eq:systemintro}), the same
holds for $Q''$ and, therefore, $Q$ is continuous at $\xi_0$.  
Hence, the first equation of (\ref{eq:systemintro}) suggests to call the discontinuity
at $\xi_0$ {\it admissible} if and only if the triple $(U(\xi_0\pm), s)$ is a shock wave
for the reduced system (\ref{eq:reduced}).

\begin{definition}\label{def:admissibleradiative}
An {\sc admissible radiative shock wave} of the hyperbolic--elliptic system 
(\ref{eq:systemintro}) is a radiative shock wave $(U,Q)$ if there exists a discrete 
set ${\mathcal J}\subset\R$ for which
$U$ is $C^1$ off\footnote{Given ${\mathcal J}\subset\R$,
we say that a function $F\,:\,\R\to\R$ is {\sc $C^k$ off $\mathcal J$} if 
for $F\in C^k([a,b])$ for any interval $[a,b]$ such that ${\mathcal J}\cap (a,b)$ is
empty.} 
 ${\mathcal J}$, $Q$ is continuous on $\R$ and $C^2$ off ${\mathcal J}$
and the triple $(U(\xi_0\pm), s)$ is a shock wave for (\ref{eq:reduced})
for any $\xi_0\in{\mathcal J}$.
\end{definition}

Later on, it will be shown that the admissible radiative shock wave we deal with have at
most one discontinuity point.
\medskip

Now we are in position to state the main problem we deal with:

\begin{quote}
{\sc Problem.}
{\sl Given a triple $u_\pm\in\R^n$, $s\in\R$, defining a shock wave for the reduced system 
(\ref{eq:reduced}), does there exist a corresponding admissible radiative shock wave
for the hyperbolic--elliptic system (\ref{eq:systemintro})?}
\end{quote}

The analogous problem has been affirmatively solved for other regularizations
of (\ref{eq:reduced}) (see \cite{MP85} and \cite{YZ00} for the viscous and the
relaxation approximation, respectively).

In the context of radiative gas model, this problem has been addressed
for the scalar case \cite{ST92,KN99} and for specific systems \cite{LCG06}.
\medskip

The result in \cite{ST92} concerns the existence of radiative shocks for
convex fluxes $f$.
Elementary computations shows that the parameters $\varepsilon$ and
$m$ in the cited paper are related with $L, G, R$ by the relations
$\varepsilon=LG/R$ and $m=\sqrt R/LG$.
Hence the existence result in \cite{ST92} can be stated as follows.

\begin{theorem} \textbf{\cite{ST92}}
Assume $f''>0$. If
 \begin{equation}\label{STcond}
  \frac{R}{L^2\,G^2}\,
   \sup_{u\in[u_+,u_-]} 
   \Bigl\{f''(u)\Bigl(f(u_-)+s(u-u_-)-f(u)\Bigr)\Bigr\}\leq\frac14,
 \end{equation}
then there exist a continuous radiative shock.
\end{theorem}

Condition (\ref{STcond}) is satisfied in the case of small shocks, i.e. if
$|u_--u_+|$ is small.  
Also, by a hyperbolic rescaling $\partial_x\mapsto L^{-1}\partial_x$,
$\partial_t\mapsto L^{-1}\partial_t$, the system (\ref{eq:systemintro}) becomes
\begin{equation*}
    \begin{cases}
	u_{t}+ f(u)_{x} +q_{x}=0 & \\
    \displaystyle{
	-m^2\,\nu^2\,L^{-2}\,q_{xx} + q +\nu\,u_{x} =0,} &
    \end{cases}
\end{equation*}
where $m^2=R/L^2G^2$ and $\nu=G\,L/R$.
As $m\to 0$, formally, we get a scalar viscous conservation law.
Hence condition (\ref{STcond}) can be read as a measure of the 
strength of the viscosity term, needed for smoothness.

In \cite{KN99}, it has been considered the scalar case with the 
choice of a Burgers' like flux function $f(u)=\frac12\,u^2$ and,
without loss of generality, $L=R=G=1$.

\begin{theorem} \textbf{\cite{KN99}} 
Let $u_\pm$ be such that $u_+<u_-$ and set $s:=(u_++u_-)/2$.
Then there exists a (unique up to shift) admissible radiative shock wave $(U,Q)$.
\end{theorem}

The authors can also determine the smoothness of the profiles $(U,Q)$ 
as the size $|u_+-u_-|$ varies, showing that regularity improves as
the amplitude of the shock decreases:\\ 
-- if $|u_+-u_-|>\sqrt 2$, then $U\in C^0(\R\setminus\{\xi_0\})$ 
for some $\xi_0\in\R$ and $Q$ is Lipschitz continuous;\\
-- if  $|u_+-u_-|<2\sqrt{2n}/(n+1)$ for some $n\in\N$, 
then $U\in C^{n}(\R)$ and $Q\in C^{n+1}(\R)$.
\medskip

In the case of systems (i.e. $u\in\R^n$), the literature on the 
subject restricts to the very
recent article \cite{LCG06}.
There, the authors consider the problem of existence of  \underline{smooth}
radiative shocks for a specific model describing the gases not in thermodynamical 
equilibrium with radiations.
In their case, the coupled elliptic equation is nonlinear with respect to the 
temperature.
\medskip

In the present article, we deal with general systems (\ref{eq:systemintro})
and consider admissible radiative shock waves, hence possibly discontinuous.
We are able to prove that, for small amplitude shock waves of (\ref{eq:reduced}), 
there always exists a radiative shock profiles.
 
\begin{theorem}\label{theo:mainintro}
Let $u_-\in\R^n$ be such that the $k-$th characteristic field of (\ref{eq:reduced})
is genuinely nonlinear at $u_-$, that is $\nabla \lambda_{k}(u_-) \cdot r_{k}(u_-)\neq 0$.
Assume that 
\begin{equation}
   (\ell_{k}(u_{-})\cdot L)\,(G\cdot r_{k}(u_{-})) >0.
\end{equation}
Then there exists a sufficiently small neighborhood $\mathcal U$ of $u_-$ such that
for any $u_+\in{\mathcal U}$, $s\in\R$ such that the triple $(u_{-},u_{+};s)$ 
defines
a shock wave for (\ref{eq:reduced}), there exists a (unique up to shift) admissible
radiative shock wave for (\ref{eq:systemintro}).
\end{theorem}

The proof is essentially divided in two steps. 
The first one is to reduce the system case to the scalar case. 
This (surprising!) possibility is essentially due to the fact that the diffusion
matrix $L\otimes G$ is rank--one. 
The second step is to generalize the existing proof for Burgers' like flux
to general strictly convex fluxes and to show the existence of an heteroclinic orbit 
connecting the asymptotic states. 
The analysis is complicated by the fact that the differential equation for the
profile is not in normal form and that discontinuity of $U$ can arise. 
In dealing with this problem, it turns to be useful to work with the integrated
variable $z$ such that $q=-z_x$, hence to deal with the equations
\begin{equation}
    \begin{cases}
	u_{t}+ f(u)_{x}=L\,z_{xx} & \\
	-z_{xx} + R\,z=G\cdot u, &
    \end{cases}
    \label{eq:integrated}
\end{equation}
instead of the original system (\ref{eq:systemintro}). 

\begin{remark}\label{rem:smalness}
The smallness assumption on the shock $(u_{-},u_{+};s)$ 
is needed only in the first step of the proof, that is in the 
reduction procedure from system to scalar case, 
which is carried out in terms of Implicit Function Theorem. 
Since we are able to prove the existence of (possibly discontinuous) radiative shocks 
for general strictly convex scalar models and general large admissible shocks,
our existence result includes in principle the case of discontinuous 
radiative shocks for (\ref{eq:systemintro}).
\end{remark}
\medskip

By applying a similar strategy to the one used in \cite{KN99}, we can also get
analogous result on smoothness of the profile $(U,Q)$. 

\begin{theorem}\label{theo:reg}    
There exists a sequence $\{\e_{n}\}$, $\e_n\to 0$ as $n\to\infty$,
such that the profile $u$ is $C^{n+1}$ whenever $|u_+-u_-|<\e_{n}$.
\end{theorem}

For completeness, once the convex case has been treated in full generality, 
we consider the case of nonconvex fluxes, again showing the existence and 
regularity of admissible radiative shock profiles. 
The number of possible jumps of the profiles and the construction of the 
profiles themselves are strictly related with the number of inflection points
of the flux function (see Section \ref{subsec:nonconvex}).
\medskip

The paper is organized as follows. 
Section \ref{sec:ex} is devoted to prove the existence result in the 
scalar case for strictly convex flux functions $f$.
Gaining of regularity of the profile as the amplitude decrease is considered 
in Section \ref{sec:reg}, again in the strictly convex case.
The scalar case with a general (nonconvex) flux function $f$ is treated 
in Section \ref{subsec:nonconvex}.
Finally, in Section \ref{sec:system}, we show how to reduce the problem from the 
system case to the scalar one.

\section{Existence of the profile in the convex case}\label{sec:ex}

We start our analysis with the  case of  a scalar conservation law, with convex flux, 
coupled with the linear elliptic equation describing the radiating effects. 
Specifically, we discuss the existence of a travelling wave profile 
for the $2\times 2$ hyperbolic--elliptic system
\begin{equation*}
    \begin{cases}
	u_{t}+ \tilde f(u)_{x} +L\,q_{x}=0 & \\
	-q_{xx} + R\,q +G\,u_{x} =0, &
    \end{cases}
\end{equation*}
where $x\in\R$, $t>0$, $u$ and $q$ are scalar functions, the flux function $\tilde f$ 
is strictly convex, $R>0$ and $L,G$ are constants such that $LG>0$.
As reported in the Introduction, previous results on the same problem are
contained in \cite{ST92} (general convex fluxes, smooth radiative shocks)
and in \cite{KN99} (Burgers' flux, general radiative shocks).

Applying the rescaling 
 $$
  \partial_t\mapsto L\,G\,\,\partial_t,\qquad 
  \partial_x\mapsto \sqrt R\,\partial_x,\qquad 
  q\mapsto G\,q/\sqrt{R},
 $$  
and setting $f(s):=\sqrt R\,\tilde f(s)/L\,G$, we get the (adimensionalized) version
\begin{equation}
    \begin{cases}
	u_{t}+ f(u)_{x} +q_{x}=0 & \\
	-q_{xx} + q +u_{x} =0, &
    \end{cases}
    \label{eq:scalar}
\end{equation}
where $f$ is a strictly convex function.

We stress that, thanks to the discussion  of Section \ref{sec:system}, the results of the 
present section will give the existence of a $k$--travelling wave solutions for system 
(\ref{eq:systemintro}) with genuinely nonlinear $k$ characteristic field.
The study of the general scalar case is left to  Section \ref{subsec:nonconvex} and the  
results for system (\ref{eq:systemintro}) without GNL assumptions is again 
a consequence of the reduction arguments of Section \ref{sec:system}.  

\begin{remark}\label{rem:hypconvexity}
    The existence of the profile does not need the strict convexity of  the flux function $f$. 
    Indeed, our proof applies also to the case of a function $f$ such that
    $f(u)-cu$ is monotone on two intervals $(u_{+},u_{*})$ and $(u_{*},u_{-})$.
\end{remark}

Let us consider a 
solution to (\ref{eq:scalar}) of the form $(u,q) = (u(x-st),q(x-st))$ such that
\begin{equation*}
    u(\pm\infty) = u_{\pm},\quad s = 
    \frac{f(u_{+})-f(u_{-})}{u_{+}-u_{-}},\ \quad u_{+}<u_{-}.
\end{equation*}
We introduce the variable $z$ as the 
opposite of the antiderivative of 
$q$, that is $z_{x} := -q$, with $z(\pm \infty) = z_{\pm} =  u_{\pm}$.
Hence, after an integration of the second equation, system (\ref{eq:scalar}) 
rewrites as
\begin{equation*}
    \begin{cases}
	u_{t}+ f(u)_{x} -z_{xx}=0 & \\
	-z_{xx} + z -u =0.
    \end{cases}
\end{equation*}
Thus, the equations for the profiles $u(x-st)$ and $z(x-st)$ take the form
\begin{equation*}
    \begin{cases}
	-su' + f(u)' - z'' = 0 & \\
	-z'' + z-u=0,
    \end{cases}	
\end{equation*}
that is, after integration of the first equation
\begin{equation}
    \begin{cases}
	-s(u-u_{\pm}) + f(u) - f(u_{\pm}) = z'  & \\
	u=z-z''.
    \end{cases}
    \label{eq:profiles}
\end{equation}
At this point, we rewrite (\ref{eq:profiles}) as the following second 
order equation for $z$
\begin{equation}
    z'=F(z-z'';s),
    \label{eq:zrewritten}
\end{equation}
where $F(u;s) = f(u) - f(u_{\pm}) -s(u-u_{\pm}) $. We shall prove the 
existence of a profile for $z$ solution of (\ref{eq:zrewritten}) 
between the states $z_{-}=u_{-}$ and 
$z_{+}=u_{+}$, $z_{-}>z_{+}$, which will give the existence of our 
profile for $u = z-z''$.
To this end, let us note that, thanks to the strict convexity of $f$, 
the function $F(\cdot;s)$ is strictly decreasing in  an interval 
$[z_{+},z_{*}]$, and strictly increasing in $[z_{*}, z_{-}]$, where 
$F(z_{\pm};s) = 0$ and $F(z_{*};s)=-m<0$. Hence, $F(\cdot;s)$ is 
invertible in the aforementioned intervals and 
we
denote with $h_{\pm}$ the corresponding inverse functions.
In the next proposition, we analyze the 
behavior of the two ordinary differential equations $z''= z-h_{\pm}(z')$ 
in the corresponding intervals of existence.

\begin{proposition}\label{prop:maximalminus}
Let us  denote with $z_+=z_{+}(x)$ the (unique up to a shift) maximal solution  of 
    \begin{equation}
	z'' =z-h_{+}(z'),
        \label{eq:minus}
    \end{equation}
with $z_{+}(+\infty) = z_{+}$ and $z_{+}'(+\infty) = 0$. 
Then $z_{+}$ is monotone decreasing,  $z_{+}'$  is monotone increasing and 
moreover $z_{+}$ is not globally defined, that is, there exists a point, assumed to be 
$0$ (thanks to translation invariance),such that
    \begin{equation}
	z_{+}(0)- z_{+}''(0) = z_{*},\quad z_{+}'(0) = - m.
        \label{eq:extrememinus}
    \end{equation}
\end{proposition}

\begin{proof}
   We rewrite (\ref{eq:minus}) as a system of first order equations 
   as follows
   \begin{equation}
       X' = H_{+}(X),
       \label{eq:minusrewrit}
   \end{equation}
   where 
   \begin{equation*}
       X = 
       \begin{pmatrix}
           z  \\
           z'
       \end{pmatrix},\quad H_{+}(X) = 
       \begin{pmatrix}
           z'  \\
           z- h_{+}(z')
       \end{pmatrix}.
   \end{equation*}
   Then 
   \begin{equation*}
        \nabla H_{+}
          \begin{pmatrix}
              z_{+}  \\
              0
          \end{pmatrix} = 
          \begin{pmatrix}
              0 & 1  \\
              1 &\displaystyle{ -\frac{1}{F'(z_{+};s)}}
          \end{pmatrix}
   \end{equation*}
and $(z_{+},0)$ is a saddle point. 
We are interested in the stable manifold at this point and, since we want $z$ to be 
decreasing, let us follow the trajectory which exits from $(z_{+},0)$ in the lower half plane 
of the states space. 
Then we claim that $z$ is monotone decreasing and $z'$ is monotone increasing. 
Indeed, if by contradiction $z$ attains a local maximum, say at $x=x_0$, then
$z'(x_0)=0$ and $0 \geq z'' = z-h_{+}(0) = z - z_{+} $ at $x=x_0$, which is impossible. 
Thus $z$ is monotone decreasing and $z'< 0$.  
Now, if by contradiction $z'$ attains a local minimum at $x=x_0$, then the 
trajectory  $z'=\varphi(z)$ in the $(z,z')$ plane must attain a local minimum at $x=x_0$, 
that is $\varphi'(z) =0$ and $\varphi''(z) \geq 0$. 
Thus, at $x=x_0$,
\begin{equation*}
 0 = \varphi'(z) = \frac{z''}{z'} = \frac{z-h_{+}(z')}{z'}
\end{equation*}
and
 \begin{align*}
\varphi''(z) &= \frac{d}{dz} \left (\frac{z - h_{+}(z')}{z'}\right ) =
\frac{(z'-h_{+}'(z')z'')z' - (z-h_{+}(z'))z''}{(z')^{3}}= \frac{1}{z'}<0,
\end{align*}
which is impossible. 
Hence  $z'$ is monotone increasing and clearly $h_{+}(z')\in [z_{+},z_{*}]$: 
Therefore $z'' = z +O(1)$, which implies the solution does not blow up in finite time.

Thus, we are left with the two following possibilities:
	\begin{itemize}
	    \item[(i)] the solution is defined for any $x\in\R$;
	
	    \item[(ii)] the solution reach the boundary of the domain of 
	    definition of the differential equation in a finite point, which 
	    can be assumed to be $0$ up to a space translation, namely, 
	    (\ref{eq:extrememinus}) holds.
	\end{itemize}
	Hence, the proof is complete if we can exclude the case (i). 
	To this 
	end, let us assume  that (i) holds. Thus in particular $z'$  
	is bounded for 
	any $x\in\R$, because $z'\in (-m,0)$. 
	Since $z$ is 
	monotone decreasing, then $z\to l\in(z_{+},+\infty]$ as $x\to-\infty$. 
	If $l<+\infty$, since $z'$ is increasing, then $z'\to 0$ as 
	$x\to-\infty$ 
	and therefore $z'' = z-h_{+}(z')\to l - h_{+}(0) = l -z_{+} = 0$, 
	which is impossible. On the other hand, let us assume 
	$l=+\infty$. Since $z'$ is globally bounded, $z$ decreases at 
	most linearly, which is incompatible with $z'' = z +O(1)$, 
	because the latter implies an exponential rate.
        Therefore, the solution cannot be defined 
	for any $x\in\R$ and the proof is complete.
\end{proof}

With the same kind of arguments, it is possible to analyze the behavior of the 
maximal solution $z_{-}(x)$ of the equation $z''  =z-h_{-}(z')$. 
Hence, the following proposition holds.

\begin{proposition}\label{prop:maximalplus}
Let us  denote with $z_{-}(x)$ the (unique, up to a space shift) maximal solution of 
    \begin{equation}
	z'' =z-h_{-}(z'),
	\label{eq:plus}
    \end{equation}
 with $z_{-}(-\infty) = z_{-}$ and $z_{-}'(-\infty) = 0$. 
Then $z_{-}$ and $z_{-}'$ are monotone decreasing and moreover $z_{-}$ is 
not globally defined, that is, there exists a point, which we can assume to be $0$
(thanks to translation invariance), such that
 \begin{equation*}
	z_{-}(0) - z_{-}''(0)   = z_{*},\quad z_{-}'(0) = -m.
 \end{equation*}
\end{proposition}

With the aid of Propositions \ref{prop:maximalminus} and \ref{prop:maximalplus}, 
we shall build a $C^{1}$ trajectory joining $z_{-}$ and $z_{+}$ by finding a point of 
intersection for the orbits of the maximal solutions $z_{-}(x)$ and $z_{+}(x)$
in the state space $(z,z')$. 
The monotonicity of $z_{\pm}(x)$ and $z_{\pm}'(x)$ shall guarantee such an intersection 
is unique and therefore the resulting $C^{1}$ trajectory from $z_{-}$ to $z_{+}$ is unique, 
up to a space translation. 
The existence of the aforementioned intersection is a straightforward 
consequence of the following lemma.

\begin{lemma}\label{lem:intersection}
    Let us  denote with $z_{+}(x)$ and $z_{-}(x)$ the maximal solutions of 
    (\ref{eq:minus}) and (\ref{eq:plus}) respectively. Then
    \begin{equation}
        z_{+}(0)\geq z_{*} \geq z_{-}(0). 
        \label{eq:intersection}
    \end{equation}
\end{lemma}

\begin{proof}
    We prove only the inequality on the left of 
    (\ref{eq:intersection}), the right one to be proved similarly.
    Let us denote with $\bar y$ the heteroclinic orbit of
    \begin{equation*}
       \begin{cases}
	   y' = F(y;s) &   \\
           y(\pm\infty) = z_{\pm}. & 
       \end{cases}
    \end{equation*}
    Our aim is to compare this solution with the solution $z_{+}(x)$ in the phase space. 
    To this end, as before we denote with $\varphi=\varphi(z)$ the function whose graph is the 
    trajectory of $(z_{+}(x), z_{+}'(x))$ in the $(z,z')$ plane 
    and with $\psi=\psi(z) = F(z;s)$ the function whose graph is the trajectory of 
    $(\bar y(x), \bar y'(x))$ in the $(z,z')$ plane. 
    Then
    \begin{equation*}
        \psi'(z_{+}) = F'(z_{+};s)<0
    \end{equation*}
    and
    $\varphi'(z_{+})= \lambda_{1}$ is the negative eigenvalue of 
    $\nabla H_{+}(z_{+},0)$, with $H_{+}$  defined in the proof of Proposition 
    \ref{prop:maximalminus}, namely, the negative root of 
    \begin{equation*}
	P(\lambda) = \lambda^{2} + \frac{\lambda}{F'(z_{+};s)} -1.
    \end{equation*}
    Since $P(F'(z_{+};s)) = (F'(z_{+};s))^{2}>0$, then 
    $F'(z_{+};s)<\lambda_{1}$ and therefore the trajectory $\psi(z)$ 
    leaves the $z$-axis in $z_{+}$ below the trajectory $\varphi(z)$. 
    
    If $\varphi(z)$ intersect $\psi(z)$ in a point of the $(z,z')$ 
    plane, then  in that point we have $\varphi(z) = \psi(z)$ and 
    $\varphi'(z) \leq \psi'(z)$. Hence, there exist a point $\bar x$ 
    such that  $\bar y(\bar x)=z_{+}(\bar x)$ and 
      $\bar y'(\bar x)=z'_{+}(\bar 
      x)<0$. Moreover, since 
      \begin{equation*}
	  \varphi'(z) = \frac{z_{+}''}{z_{+}'} 
      \end{equation*}
      and
    \begin{equation*}
	\psi'(z) = \frac{d}{dz} \bar y' = \frac{\bar y''}{\bar y'},
    \end{equation*}
      in that point we  also have  $\bar y'' \leq z_{+}'' = z_{+} - 
      h_{+}(z'_{+}) = \bar y - h_{+}(\bar y') =0$. On the other hand, 
      $0\geq \bar y'' = F'(\bar y;s)\bar y'= F'(\bar y;s)F(\bar y;s)$, that is  
      $F'(\bar y;s)\geq 0$ because $F(\bar y;s)<0$, namely $\bar y = 
      z_{+} \geq z_{*}$ in that point. 
      
      Conversely, let us assume that $\varphi(z)$ remains above $\psi(z)$. 
      From Proposition 
      \ref{prop:maximalminus} we know that  $z'_{+}(0) = -m$, which is 
      the minimum of $\psi(z)$, attained at $z_{*}$. Therefore 
       we must have $z_{+}(0) = z_{*}$ and the proof is complete.
\end{proof}

With the aid of Lemma \ref{lem:intersection}, we are able to prove the 
main result of this section, which is contained in the following 
theorem.

\begin{theorem}\label{theo:existence}
There exists a (unique up to space translations) $C^{1}$ profile $z$ 
 with $z(\pm\infty) = z_{\pm}=u_\pm$ such that the function $z(x-st)$ is
solution of (\ref{eq:zrewritten}), where $s$ is given by the Rankine--Hugoniot condition.
The solution $z$ is $C^{2}$ away from a single point, where $z''$ has at most a jump discontinuity.
   
Moreover, there exist a (unique up to space translations) profile $u$ with 
$u(\pm\infty) = u_{\pm}$ such that the function  $u(x-st)$ is 
solution of (\ref{eq:profiles}), where $s$ is given by the Rankine--Hugoniot condition.  
This profile is continuous away from a single point, where it has at most a jump 
discontinuity which verifies the Rankine--Hugoniot and the admissibility conditions of 
the scalar conservation law $u_{t}+f(u)_{x}=0$.
\end{theorem}

\begin{proof}
    Let us observe that Lemma \ref{lem:intersection}
    implies that there exists a point in the 
    $(z,z')$ plane where the graphs of $z_{+}$ and $z_{-}$ intersects, 
    namely $z_{-}=z_{+}$ and $z_{-}'=z_{+}'$ in that point. 
    Moreover, due to the monotonicity of these graphs, which comes from 
    the proofs of Proposition \ref{prop:maximalminus} and Proposition 
    \ref{prop:maximalplus}, this intersection is indeed unique.  
    Hence with  appropriate space translations, we can find a point $\bar x$ 
    such that $z_{-}(\bar x) = z_{+}(\bar x) = \bar z$ and $z'_{-}(\bar x) = 
    z'_{+}(\bar x)= \tilde z$. 
   Set  
    \begin{equation*}
        z(x) := 
        \begin{cases}
            z_{-}(x) &  x\leq \bar x \\
            z_{+}(x) & x\geq \bar x.
        \end{cases}
    \end{equation*}
    Thus, this profile is the unique (up to space translations) $C^{1}$ solution of
    (\ref{eq:zrewritten}) and it verifies $z(\pm\infty) = z_{\pm}$. 
Moreover, $z(x)$ is $C^{2}$ in the intervals $(-\infty,\bar x]$ and $[\bar x,+\infty)$ and 
finally $z''(\bar x -0) = z''_{-}(\bar x) = z_{-}(\bar x) - h_{-}(z_{-}'(\bar x)) =\bar z - h_{-}(\tilde z)$ 
and $z''(\bar x +0) = z''_{+}(\bar x) = z_{+}(\bar x) - h_{+}( z_{+}'(\bar x)) =\bar z - h_{+}(\tilde z)$.
   
   The regularity of $u=u(x-st)$ is a direct consequence of 
   the first part of the theorem and of the relation $u=z-z''$. 
   Moreover, in the case of a discontinuity in $u$, namely $u(\bar x 
   -0)\neq u(\bar x +0)$,  $(u(\bar x -0), u(\bar x +0); s)$
   verifies 
      the Rankine--Hugoniot condition 
      for the strictly 
	 convex conservation law $u_{t}+f(u)_{x}=0$. Indeed, $u(\bar x -0) = 
   h_{-}(\tilde z)$, $u(\bar x +0) = h_{+}(\tilde z)$ and a direct 
   calculation shows
   \begin{equation*}
       \frac{f(h_{+}(\tilde z))-f(h_{-}(\tilde z))}{ h_{+}(\tilde z)- h_{-}(\tilde z)} = 
        s.
   \end{equation*}
   Finally, $u(\bar x -0) = h_{-}(\tilde z) > h_{+}(\tilde z)= 
   u(\bar x +0)$, that is, this shock is admissible and the proof is 
   complete.
\end{proof}

\begin{remark}\label{rem:uniqueness}
    The above theorem contains an uniqueness result, among the class 
    of radiative shocks given in Definition \ref{def:radiative}. An
    uniqueness result in a wider class of solutions  is contained in 
    Theorem \ref{theo:uproperties}. We postpone this result because 
    it is  valid for general 
    flux functions, disregarding convexity properties.
\end{remark}    

\begin{remark}\label{rem:regularity}
   From the last lines of the above proof, it is clear that the 
   profile $u$ can have a jump discontinuity in a point, as it was 
   already proved in \cite{KN99} for sufficiently large shocks in the 
   case of the Burgers' equation. Moreover, in that paper it is proved 
   also that, below an explicit threshold, the profile is continuous 
   and it smoothes out  as the strength of the shock decreases.
   In our case, the profile in $u$ is continuous if $z''(\bar x -0) = 
   z''(\bar x +0)$, namely when $h_{-}(\tilde z) =h_{+}(\tilde z) = 
   z_{*}$, that is $-m =\tilde z =z_{\pm}'(\bar x)$, $z_{\pm}(\bar x) = 
   z_{*}$ and $z_{\pm}''(\bar x) =0$. 
   This property and the further regularity of the profile for sufficiently small shocks
   is proved in Section \ref{sec:reg} below.
\end{remark} 

\section{Regularity of the profile in the convex case}\label{sec:reg}
In Section \ref{sec:ex}, we proved the existence of  the travelling wave 
$(u,q) = (u(x-st),q(x-st))$,
\begin{equation*}
    u(\pm\infty) = u_{\pm},\quad s =
    \frac{f(u_{+})-f(u_{-})}{u_{+}-u_{-}},
\end{equation*}
for (\ref{eq:scalar}), when the flux $f(u)$ is smooth and strictly
convex; in the present section we focus our attention to its regularity.

Let us start by observing that this issue is related to
the smallness of the shock we are dealing with \cite{KN99}.
Here we shall prove that it is possible to recast such a property in
terms of smallness of the term $-q_{xx}$ in (\ref{eq:scalar}), namely 
when this system increases its diffusive
nature, being closer to its parabolic equilibrium. This will be made 
in terms of a diffusive scaling with respect to the shock strength 
$\e = |u_{+}-u_{-}|$
and, following
the ideas of \cite{KN99},  by analyzing the system for $(u,v=u')$. In 
that way, we will obtain the same kind of results of 
\cite{KN99} for the special case $f(u)= \frac{1}{2}u^{2}$. 
Let us finally observe that  similar phenomena arise in the
discussions of  Section \ref{subsec:nonconvex} to prove the
existence and regularity of travelling wave profiles of (\ref{eq:nonconvex}) with  
flux functions $f$ with change of convexity. In that case, the small
parameter  $\e$ is already present in the model and it is not 
connected with the smallness of the shock. 
However, existence and 
regularity of the profile will require once again the smallness of 
that parameter,
namely, as before, when the hyperbolic--elliptic model is close enough
to its diffusive underlying dynamic.  

From (\ref{eq:scalar}) and after integration with respect to $\xi$ of 
the first line, the equations for the profiles $u$ and $q$ are given by
\begin{equation}
    \begin{cases}
	-s(u-u_{\pm}) + f(u) - f(u_{\pm}) + q = 0 & \\
	u'=q''-q,
    \end{cases}
    \label{eq:uqprofile}
\end{equation}
because $q(\pm\infty) =0$. Thus, using $q = - [f(u) - f(u_{\pm})
-s(u-u_{\pm})] = -F(u;s)$ in (\ref{eq:uqprofile})$_{2}$, we end
up with
\begin{equation*}
    u' =- F(u;s)'' + F(u;s) = -(f'(u) -s)u''
    -f''(u)(u')^{2} + F(u;s).
\end{equation*}
Hence, we obtain the following system in the state space $(u,v=u')$
\begin{equation}
    \begin{cases}
	u'= v & \\
	(f'(u) -s)v'=-f''(u)v^{2}-v+F(u;s).
    \end{cases}
    \label{eq:uvprofile}
\end{equation}
It is worth to observe that system (\ref{eq:uvprofile}) is singular
where $(f'(u) -s)=0$  and, thanks to the strict convexity of
$f$, the latter occurs at a unique $\bar u\in (u_{+},u_{-})$.
At this stage, driven by
\cite{LM03}, we scale both the independent and the dependent variable 
as follows:
\begin{align*}
    & \widetilde{u}  = \frac{u-u_{+}}{\e} & & \widetilde{v} =
    \frac{v}{\e^{2}} & & \tilde\xi =\e\xi.
\end{align*}
Then, dropping the tildas, system (\ref{eq:uvprofile}) becomes
\begin{equation}
    \begin{cases}
	u'= v & \\
	f_{\e}'(u)v'=\frac{1}{\e^{2}}\left (-\e^{2}f_{\e}''(u)v^{2}-v+f_{\e}(u)\right
	),
    \end{cases}
    \label{eq:uvscaled}
\end{equation}
where $f_{\e}(u):= \frac{1}{\e^{2}}F(u_{+}+\e u;s)$ and $u\in [0,1]$. 
Once again, system (\ref{eq:uvscaled}) is singular at $u=\bue$, for a unique $\bue\in (0,1)$.
Finally, as in \cite{KN99}, we remove the singularity of (\ref{eq:uvscaled}) by using once 
again a new independent variable  $\eta$ defined by
\begin{equation*}
    \xi = \int_{\eta}^{\infty} f_{\e}'(u(\zeta))d\zeta.
\end{equation*}
 The resulting system reads as follows
\begin{equation}
    \begin{cases}
	u'= f_{\e}'(u)v & \\
	v'=\frac{1}{\e^{2}}\left (-\e^{2}f_{\e}''(u)v^{2}-v+f_{\e}(u)\right
	).
    \end{cases}
    \label{eq:uvscaled2}
\end{equation}
We notice that (\ref{eq:uvscaled2}) admits $(0,0)$ and $(1,0)$ as  equilibrium
points, corresponding to the two equilibrium points  $(u_{\pm},0)$ of 
the original system.
\begin{remark}\label{rem:vari}
(i) For any $n\geq 2$, if the original flux function $f$ is $C^{n}$, then $f_{\e}\to
    \frac{1}{2}f''(u_{+})u(u-1)$ in $C^{n}([0,1])$, as $\e\downarrow 0$.\\ 
(ii) The last change of independent variable gives a
	reparametrization of the orbit of (\ref{eq:uvscaled}) in the
	two regions $[0,\bue)$ and $(\bue,1]$. Let us first observe that 
	$f_{\e}'(u)<0$ for $0 \leq u<\bue$ and $f_{\e}'(u)>0$ for $\bue < u \leq 1$.
	Hence,
	we obtain continuity of the orbit $u$ of (\ref{eq:uvscaled}), provided we prove the 
	two  previous reparametrized orbits verify $u(-\infty) =
	0$, $u(+\infty) = \bue$, and  
	   $u(+\infty) =
	\bue$, $u(-\infty) = 1$ respectively. We shall prove the
	former property, the latter being similar.\\
(iii)  As for \cite{KN99}, further regularity for the profiles
	solutions of (\ref{eq:uvscaled2}) implies further regularity of the original
	profile, thanks to their exponential
	decay toward the asymptotic states at $\pm\infty$.
	\end{remark}
For the sake of clarity, we state first the result concerning  the continuity and the $C^{1}$ 
regularity for the profile $u$ solution of (\ref{eq:uvscaled2}) and then the one for the 
$C^{2}$ regularity. 
The general case, including the former ones, is then stated in the last Proposition. 
Here below, we shall assume $f$ to be smooth.

\begin{proposition}\label{theo:conteC1}
    There exist two values $\bar\e$, $\e_{0}>0$ such that, for $\e <
    \min\{\bar\e, \e_{0}\}$, the orbits of (\ref{eq:uvscaled2}) which
    pass through the equilibrium points $(0,0)$ and $(1,0)$ meet at the
    equilibrium
    point $(\bue, \bve)$, where
    \begin{equation} 
	\bve = \frac{-1 + \sqrt{1+
	4\e^{2}f''_{\e}(\bue)f_{\e}(\bue)}}{2\e^{2}f''_{\e}(\bue)}.
	\label{eq:newv2}
    \end{equation}
    In particular, the orbit $u$ is $C^{1}$.
\end{proposition}

\begin{proof}
    We start by studying the ``new'' equilibrium points of
    (\ref{eq:uvscaled2}) introduced by the last change of variable, besides the 
    aforementioned points $(0,0)$ and $(1,0)$. Hence we must satisfy
    the relations 
    \begin{equation*}
        \begin{cases}
            f'_{\e}(u)=0 &   \\
	    \e^{2}f_{\e}''(u)v^{2}+v-f_{\e}(u)=0,& 
        \end{cases}
    \end{equation*}
    that is $u=\bue$ and 
    \begin{equation*} 
		\bar v^{\e}_{1,2} = \frac{-1 \mp \sqrt{1-
		4\e^{2}f''_{\e}(\bue)|f_{\e}(\bue)|}}{2\e^{2}f''_{\e}(\bue)},
	    \end{equation*}
	    provided $\e\leq \e_{0}$,
    thanks to Remark \ref{rem:vari}--(i). Moreover, using again Remark
    \ref{rem:vari}--(i), there exists an $\bar\e>0$ such that for any
    $\e\leq \bar\e$  the whole curves
    \begin{equation} 
	    v^{\e}_{1,2}(u) = \frac{-1 \mp \sqrt{1-
	    4\e^{2}f''_{\e}(u)|f_{\e}(u)|}}{2\e^{2}f''_{\e}(u)}
	    \label{eq:curves}
	\end{equation}
with $u\in [0,1]$ are real. We shall now study the nature of these
points to prove the first assertion of the theorem. The Jacobian associated to system
(\ref{eq:uvscaled2}) is given by
\begin{equation*}
    J(u,v) = \begin{pmatrix}
        f''_{\e}(u)v & f'_{\e}(u)  \\ \\
        -f'''_{\e}(u) +\frac{1}{\e^{2}}f'_{\e}(u) & -2f''_{\e}(u)v -
	\frac{1}{\e^{2}}
    \end{pmatrix}
\end{equation*}
Evaluating $J(0,0)$, $J(1,0)$ and $J(\bue, \bve)$ we conclude $(0,0)$ 
and $(1,0)$ are saddle points, while $(\bue,\bve)$ is a sink. 
Moreover, the curve $v= v^{\e}_{2}(u)$ given in (\ref{eq:curves})
passes through the point $(0,0)$ and
\begin{equation*}
    \left.\frac{d}{du}v^{\e}_{2}(u)\right|_{u=0} = f'_{\e}(0).
\end{equation*}
Let us consider the orbit exiting from $(0,0)$ at $-\infty$, namely
the one on the unstable manifold of that (saddle) point, in the
halfplane $v<0$.  Its tangent 
vector is then given by the eigenvector related to the positive
eigenvalue $\lambda_{+}$ of $J(0,0)$, that is 
\begin{equation*}
   \begin{pmatrix} 1  \\ \frac{\lambda_{+}}{f'_{\e}(0)}  \end{pmatrix} 
   = \begin{pmatrix} 1  \\ \frac{2f'_{\e}(0)}{1+ \sqrt{1+4\e^{4}f'_{\e}(0)}} \end{pmatrix}
\end{equation*}
A direct calculation shows 
\begin{equation*}
    0>\frac{\lambda_{+}}{f'_{\e}(0)}> f'_{\e}(0),
\end{equation*}
namely the trajectory leaves the equilibrium point above the
curve $v= v^{\e}_{2}(u)$. Hence, a straightforward analysis of the
dynamical system (\ref{eq:uvscaled2}) implies it must reach the
equilibrium point $(\bue, \bve)$ at $+\infty$. 
The same result holds for the orbit exiting from $(1,0)$ and therefore the first claim of 
the theorem is proved. 
Finally, since $v$ is continuous and bounded everywhere,  $u\in C^{1}$
and the proof is complete.
\end{proof}

\begin{proposition}\label{theo:C2}
    There exists a value $0<\e_{1}<\e_{0}$ such that for $\e <
    \min\{\bar\e, \e_{1}\}$, the orbit $u$ of (\ref{eq:uvscaled2}) is
    $C^{2}$.
\end{proposition}

\begin{proof}
    We evaluate the Jacobian $J(u,v)$ at the equilibrium point $(\bue,\bve)$ to obtain
    \begin{equation*}
	J(\bue,\bve) = \begin{pmatrix}
	    f''_{\e}(\bue)\bve & 0  \\ \\
	    -f'''_{\e}(\bue)\left (\bve\right)^{2}   & -2f''_{\e}(\bue)\bve -
	    \frac{1}{\e^{2}}
	\end{pmatrix}
    \end{equation*}
    which admits the following (negative) eigenvalues
 \begin{align*}
      &  \lambda_{1}^{\e} = f''_{\e}(\bue)\bve = \frac{-1+
      \sqrt{1-4\e^{2}f''_{\e}(\bue)|f_{\e}(\bue)|}}{2\e^{2}} \\
      & \lambda_{2}^{\e} = -2f''_{\e}(\bue)\bve -
	    \frac{1}{\e^{2}} = - \frac{
      \sqrt{1-4\e^{2}f''_{\e}(\bue)|f_{\e}(\bue)|}}{\e^{2}}.
 \end{align*}
 For $\e=\e_{0}$, $\lambda_{2}^{\e}=0$ and
 $\lambda_{1}^{\e}=-\frac{1}{2\e^{2}}<0$, while $\lambda_{2}^{\e}\to
 -\infty$ and $\lambda_{1}^{\e}\to -f''_{\e}(\bue)|f_{\e}(\bue)| <0$
 as $\e\downarrow 0$. Thus, there exists an $\e_{1}<\e_{0}$ such that 
 for $\e < \e_{1}$, $\lambda_{2}^{\e} < \lambda_{1}^{\e}<0$.
 Therefore, for any $\e < \e_{1}$, the orbit in $(\bue,\bve)$ is
 tangent to the eigenvector related to 
 $\lambda_{1}^{\e}$, that is 
 \begin{equation*}
     \begin{pmatrix}
         1  \\
         \frac{f'''_{\e}(\bue)\left (\bve\right 
	 )^{2}}{-3f''_{\e}(\bue)\bve
	- \frac{1}{\e^{2}} }
     \end{pmatrix} = 
     \begin{pmatrix}
	      1  \\
	      \frac{f'''_{\e}(\bue)\left (\bve\right
	      )^{2}}{\lambda^{\e}_{2}-\lambda_{1}^{\e} }
	  \end{pmatrix}
 \end{equation*}
 instead of the eigenvector related to $\lambda_{2}^{\e}$, namely $(0,1)^t$.
 
  Finally, 
  \begin{equation}
      \left.\frac{dv}{du}\right |_{u=\bue} = \frac{f'''_{\e}(\bue)\left (\bve\right
	      )^{2}}{\lambda^{\e}_{2}-\lambda_{1}^{\e} } =
	     \bar  w^{\e}_{1}\in\R,
	     \label{eq:w1}
  \end{equation}
  which implies the $C^{1}$ regularity of the $v$--component of the orbit 
  of (\ref{eq:uvscaled2}), 
  that is, the $C^{2}$  regularity of the  $u$--component of that orbit.
\end{proof}
To prove further regularity of the profile $(u,v)$ solution of
(\ref{eq:uvscaled2}), we follow again the ideas of \cite{KN99}. 
For $n\geq 0$, let us consider by induction  the new variable
\begin{equation}
    w_{n} = \frac{1}{(u-\bue)^{n}}\left (v - \sum_{j=0}^{n}\bar
    w_{j}^{\e}(u-\bue)^{j} \right ),
     \label{eq:newvar}
\end{equation}
that is 
\begin{equation}
    v = \sum_{j=0}^{n}\bar w_{j}^{\e}(u-\bue)^{j} + w_{n}(u-\bue)^{n} 
    : = h_{n}(u,w_{n}).
    \label{eq:change}
\end{equation}
Let us note that, with the notation introduced above, we choose $\bar w_{0}^{\e} = \bve$
and $\bar w_{1}^{\e}$   given by (\ref{eq:w1}). Hence, for $n=1$,
(\ref{eq:newvar}) reduces to
\begin{equation*}
    w_{1} = \frac{v-\bve}{u-\bue} - \bar w_{1}^{\e},
\end{equation*}
that is, $\bar w_{1}^{\e}$ is the value of the first derivative of the curve
$v=v(u)$ for $u=\bue$ and therefore the rest point $(\bue,\bve)$
becomes $(\bue,0)$ in the new pair of variables $(u,w_{1})$.
However, for $n>1$, we   assume by induction the constants  $\bar w_{j}^{\e}$, $j=1,\dots,n-1$ 
in (\ref{eq:newvar}) and (\ref{eq:change}) are given, but we determine $w_{n}^{\e}$ by 
imposing a different condition (see (\ref{eq:wj}) below) and we obtain the relation
\begin{equation}
    \bar w_{j}^{\e} = \frac{1}{j!}\left.\frac{d^{j}v}{du^{j}}\right
    |_{u=\bue},\ \hbox{for any}\ j\geq 2
    \label{eq:wjold}
\end{equation}
as a consequence.
The system for $(u,w_{n})$ is given by
\begin{equation}
    \begin{cases}
         u' =f_{\e}'(u) h_{n}(u,w_{n})&   \\ \\
        w_{n}' = \displaystyle{\frac{\theta_{n}(u,w_{n})}{\e^{2}(u-\bue)^{n}}}, & 
    \end{cases}
    \label{eq:systuwn}
\end{equation}
where
\begin{align*}
    \theta_{n}&(u,w_{n}) = -\e^{2}f''_{\e}(u)h_{n}(u,w_{n})^{2} - h_{n}(u,w_{n})
    + f_{\e}(u) \\
    &\ - \e^{2}\left ( \sum_{j=1}^{n}j\bar
    w^{\e}_{j}(u-\bue)^{j-1} + n w_{n}(u-\bue)^{n-1}\right
    )f'_{\e}(u)h_{n}(u,w_{n}).
\end{align*}
We are ready now to prove the main result of this section, namely, 
Proposition \ref{prop:reg}, (already 
proved for $n=0,1$ in Propositions \ref{theo:conteC1} and \ref{theo:C2}).

\begin{proposition}\label{prop:reg}
    For any $n\geq 0$, There exists a decreasing sequence of positive 
    values $\{\e_{n}\}_{n\geq 0}$, such that, for any $\e < 
    \min\{\bar\e, \e_{n}\}$, $v$ is a $C^{n}$ function of $u$ and 
    admits the expansion
    \begin{equation}
	v = \sum_{j=0}^{n}\bar w_{j}^{\e}(u-\bue)^{j} + 
	o((u-\bue)^{n}),
	\label{eq:exp}
    \end{equation}
    for $u\to\bue$. 
\end{proposition}

\begin{proof}
    A direct calculation shows
    \begin{equation*}
	\theta_{n}(u,w_{n}) = F_{n}(u) + G_{n}(u)w_{n} + H_{n}(u) 
	w_{n}^{2},
    \end{equation*}
    where
    \begin{align*}
	 F_{n}(u)& = -\e^{2} f''_{\e}(u) \left (  \sum_{j=0}^{n} \bar 
	 w_{j}^{\e} (u-\bue)^{j} \right)^{2} -  \sum_{j=0}^{n} \bar 
	 w_{j}^{\e} (u-\bue)^{j} + f_{\e}(u) \\
	 & \ -\e^{2} f'_{\e}(u)\left (  \sum_{j=1}^{n} j\bar 
	 w_{j}^{\e} (u-\bue)^{j-1} \right ) \left (\sum_{l=0}^{n} \bar 
	 w_{l}^{\e} (u-\bue)^{l} \right );
\end{align*}
\begin{align*}
	 G_{n}(u)& = -2\e^{2} f''_{\e}(u) \left ( \sum_{j=0}^{n} \bar 
	 w_{j}^{\e} (u-\bue)^{j+n}\right )  - (u-\bue )^{n}   \\
	 & \ -\e^{2} f'_{\e}(u)\left (  \sum_{j=1}^{n} (j+n)\bar 
	 w_{j}^{\e} (u-\bue)^{j+n-1} \right ) -\e^{2} n f'_{\e}(u)\bar 
	 w_{0}^{\e}(u-\bue)^{n-1};
	 \end{align*}
	$$
	H_{n}(u)  = -\e^{2} \left ( f''_{\e}(u) (u-\bue)^{2n} + n 
	f'_{\e}(u)(u-\bue)^{2n-1} \right ).
  $$
    In addition, for sufficiently smooth fluxes $f$, since $f'_{\e}(\bue) 
    =0$, we obtain, for $u\to\bue$,
    \begin{equation*}
	G_{n}(u) = O((u-\bue)^{n}),\quad H_{n}(u) = O((u-\bue)^{2n}).
    \end{equation*}
    Moreover,   the expression of $F_{n}(u)$ implies $F_{n}(u) = 
    F_{n-1}(u) + O((u-\bue)^{n})$ and therefore the coefficients of the 
    Taylor approximation of $F_{n}(u)$ about $\bue$ do not depend on 
    $n$, that is
    \begin{equation*}
	F_{n}(u) = \sum_{j=0}^{n}c_{j}(u-\bue)^{j} + O((u-\bue)^{n+1}),
    \end{equation*}
    with $c_{j}$ independent from $n$, for any $j$. We write down 
    these coefficients more explicitly as follows
    \begin{equation*}
	c_{j} = \alpha_{j}(\bar w_{0}^{\e},\dots, \bar w_{j-1}^{\e}) + 
	\beta_{j}(\bar w_{0}^{\e})\bar w_{j}^{\e},
    \end{equation*}
    where 
    \begin{equation*}
	\beta_{j}(\bar w_{0}^{\e}) = - \left (\e^{2} (2+j)f''_{\e}(\bue)\bar 
	w_{0}^{\e} +1 \right ).
    \end{equation*}
    Now, using again Remark \ref{rem:vari}--(i), there exists a 
    decreasing sequence  $\{\e_{j}\}_{j\geq 0}$ such that
    \begin{equation}
	\beta_{j}(\bar w_{0}^{\e})<0
        \label{eq:ineqj}
    \end{equation}
    for any 
    $\e<\e_{j}$,   and therefore we choose
      $\bar w_{j}^{\e}$, for $\e<\e_{j}$ and for any $j$, such that 
     \begin{equation}
	   c_{j} = 0\ \Leftrightarrow\ F_{n}^{(j)}(\bue)=0.
	 \label{eq:wj}
     \end{equation}
     It is worth to observe that, if we characterize the above argument, 
     and in particular the values $\e_{j}$ and (\ref{eq:wj}), 
     to the cases $n=0$ and $n=1$, we recover the previous choices 
     $\bar w_{0}^{\e} = \bve$ and $\bar w_{1}^{\e}$ given by 
     (\ref{eq:w1}). 
     
     At this point, we have determined all coefficients $\bar 
     w_{j}^{\e}$, $j\geq 0$ in (\ref{eq:exp}) and we are left to the 
     proof of the regularity of $v$ as a function of 
     $u$
     close to $u=\bue$. 
     Assume by induction the result holds for $n$ and introduce the 
     new variable $w_{n}$ defined in (\ref{eq:newvar}), with
  $\bar w_{j}^{\e}$, $j\leq n$ verifying 
      (\ref{eq:wj}). Therefore, $w_{n}\to 0$ as $u\to\bue$ by the 
      induction hypothesis and, thanks to (\ref{eq:wj}), $(\bue,0)$ is a rest point 
      for (\ref{eq:systuwn}), which correspond to the original rest point 
     $(\bue, \bve)$ for (\ref{eq:uvscaled2}).
     Since $h_{n}(\bue,w_{n})\equiv 0$ for any $w_{n}$ and in view of 
     the structure of $\theta_{n}(u,w_{n})$ showed before,
       the Jacobian of (\ref{eq:systuwn}) evaluated at this rest 
      point is given by
      \begin{equation}
          \begin{pmatrix}
              f''_{\e}(\bue) \bar w_{0}^{\e} & 0  \\ \\
	      \displaystyle{\lim_{u\to\bue} \frac{F'_{n}(u) (u-\bue)  - 
	      nF_{n}(u) }{\e^{2}(u-\bue)^{n+1}}} & 
	      \displaystyle{\lim_{u\to\bue}\frac{G_{n}(u)}{\e^{2}(u-\bue)^{n}}}
          \end{pmatrix}.
          \label{eq:jacn1}
      \end{equation}
      Since we know (\ref{eq:wj}) for $j\leq n$, then 
      \begin{equation*}
	  \lim_{u\to\bue} \frac{F'_{n}(u) (u-\bue)  - 
			nF_{n}(u) }{\e^{2}(u-\bue)^{n+1}} = s_{n}\in\R.
      \end{equation*}
      Moreover
      \begin{equation*}
	  \lim_{u\to\bue}\frac{G_{n}(u)}{\e^{2}(u-\bue)^{n}} = 
	  -\left ((2+n) f''_{\e}(\bue)\bar w_{0}^{\e} 
	  +\frac{1}{\e^{2}}\right ) = 
	  \lambda_{2}^{\e} -n\lambda_{1}^{\e}.
      \end{equation*}
      Therefore, the eigenvalues of (\ref{eq:jacn1}) are given by 
      $f''_{\e}(\bue) \bar w_{0}^{\e} = \lambda _{1}^{\e}$ and
      $\lambda_{2}^{\e} -n\lambda_{1}^{\e}$ with eigenvectors
      \begin{equation*}
          \begin{pmatrix}
              1  \\
              \displaystyle{\frac{s_{n}}{(n+1)\lambda_{1}^{\e}-\lambda_{2}^{\e}}}
          \end{pmatrix} 
            \quad\textrm{and}\quad
	  \begin{pmatrix} 0  \\ 1 \end{pmatrix}. 
      \end{equation*}
      Thus, the result for $n+1$ is obtained proceeding as in Proposition 
      \ref{theo:C2}, by imposing $\lambda_{1}^{\e} > \lambda_{2}^{\e} 
      -n\lambda_{1}^{\e}$, that is, condition (\ref{eq:ineqj}) for $j=n+1$.
      Finally, from the uniqueness of the Taylor expansion of $v$ about $\bue$, 
     relations (\ref{eq:wjold}) is also verified. 
\end{proof}

\begin{remark}\label{rem:burgers}
    In the case of a Burgers' flow, that is, the case studied in 
    \cite{KN99}, our scaling procedure will lead to a modified quadratic  flux 
    \emph{independent from} $\e$, namely $\frac{1}{2}u(u-1)$ (see 
    also \cite{LM03} for further comments on the behavior of the flux 
    with respect to this scaling). Therefore, for that particular 
    case, the thresholds $\e_{j}$ 
    can be explicitly calculated from (\ref{eq:ineqj}),   
    $\bar \e = \e_{0}$ and our results coincide with the ones 
    established in  \cite{KN99}.
\end{remark}

\section{The non convex case}\label{subsec:nonconvex}

Now, let us consider existence and regularity of travelling wave profiles when the flux 
function may change its convexity. 
As already pointed out in Section \ref{sec:reg}, in that case, we shall prove these results 
if the behavior of (\ref{eq:scalar}) is  close enough to its underlying diffusive equilibrium, 
that is, when $-q_{xx}$ is sufficiently small.
Hence, let us consider the $2\times 2$ system
\begin{equation}
    \begin{cases}
	u_{t}+ f(u)_{x} +q_{x}=0 & \\
	-\e q_{xx} + q +u_{x} =0, &
    \end{cases}
    \label{eq:nonconvex}
\end{equation}
with $0<\e \ll 1$. Let $(u_{-},u_{+};s)$ be an admissible shock for the inviscid 
conservation law $u_{t} + f(u)_{x}=0$, namely, it  satisfies 
the Rankine--Hugoniot condition
\begin{equation*}
    s = \frac{f(u_{+})-f(u_{-})}{u_{+}-u_{-}}
\end{equation*}
and the 
the strict Oleinik condition
\begin{equation}
    \frac{f(u)-f(u_{-})}{u-u_{-}}> 
    \frac{f(u_{+})-f(u_{-})}{u_{+}-u_{-}} > 
    \frac{f(u_{+})-f(u)}{u_{+}-u},
    \label{eq:ole}
\end{equation}
for any $u$ between $u_{-}$ and $u_{+}$. 
As before, let us introduce a new variable $z$ as follows $-z_{x}:=q$ 
so that the equations for $u(x-st)$ and $z(x-st)$ are given by
\begin{equation}
    \begin{cases}
	F(u;s) =z' & \\
	-\e z'' + z =u, &
    \end{cases}
    \label{eq:nonconvexprof}
\end{equation}
where, as usual, $F(u;s) = f(u)-f(u_{\pm}) - s(u-u_{\pm})$ and $z(\pm\infty) = z_{\pm} = u_{\pm}$. Without loss of generality, assume  $u_{-}>u_{+}$. 
Then (\ref{eq:ole}) becomes
\begin{equation}
    F(u;s)<0,\ \hbox{for any}\ u\in(u_{+},u_{-})
    \label{eq:ole2}
\end{equation}
and in addition it implies the Lax condition
\begin{equation*}
    f'(u_{+}) \leq s\leq f'(u_{-}).
\end{equation*}
We start by considering the non degenerate situation, that is
\begin{equation}
 \aligned
    &f'(u_{+}) < s <f'(u_{-}),\\ 
    &f''(u)\neq 0\ \hbox{for any}\ 
    u\in(u_{+},u_{-})\ \hbox{with}\ f'(u)=0.
       \endaligned
    \label{eq:laxstrict}
\end{equation}
The result in the general case will be proved by an approximation procedure at the 
end of the section \cite{Ser06}.
We shall construct a profile for the $z$--component by solving the equation
\begin{equation}
    z' = F(z-\e z'';s) 
    \label{eq:original}
\end{equation}
and then  $u$ will be given by $u=z-\e z''$. 

Since $f$ is not strictly convex, in order to invert the function $F(\cdot;s)$, we decompose 
$[z_{+},z_{-}]$ in the disjoint (up to vertices) union of subintervals $I_{1},\ldots,I_{2n}$, 
$n>1$, where  $F(\cdot;s)$ is monotone.
More precisely, thanks to (\ref{eq:ole2}) and (\ref{eq:laxstrict}), $F(\cdot;s)$ is decreasing in 
$I_{2k-1}$ and increasing in $I_{2k}$, $k=1,\ldots,n$. In the spirit of Propositions 
\ref{prop:maximalminus}  and \ref{prop:maximalplus}, we shall construct 
\emph{maximal solutions} $z_{0}$, $z^{l}_{k}$, $z^{r}_{k}$, $k=1,\ldots,n-1$ and $z_{n}$ that 
correspond respectively to $I_{1}$, $I_{2k}$, $I_{2k+1}$, 
$k=1,\ldots,n-1$ and $I_{2n}$ and  then match the corresponding graphs in 
the phase plane to obtain a global $C^{1}$ profile as in Theorem \ref{theo:existence}. 
In addition, this solution will be $C^{2}$ in the points $\{z^{*}_{2k}\} = I_{2k}\cap I_{2k+1}$, 
$k=1,\ldots,n-1$, where $F(\cdot;s)$ 
attains a local maximum, while the second derivative will have a jump 
discontinuity in  the points $\{z^{*}_{2k-1}\} =I_{2k-1}\cap I_{2k}$, 
$k=1,\ldots,n$, where $F(\cdot;s)$ 
attains a local minimum ($z^{*}_{k}$, $k=1,\ldots,2n$ are the zeros 
of $F'(\cdot;s)$ in $[z_{+},z_{-}]$, ordered from the left to the 
right). In this way, the resulting profile for $u$ 
will be regular, except for $n$ points, where it has at most a jump 
discontinuity, which is an inviscid shock satisfying (\ref{eq:ole2}).
Clearly, the the case $n=1$ correspond to the strict convex case treated before. 
Moreover, the construction of the (unique, up to a space translation) maximal 
solutions $z_{+}$ and $z_{-}$ of Propositions \ref{prop:maximalminus} 
and \ref{prop:maximalplus} can be repeated under assumptions (\ref{eq:ole2}) 
and (\ref{eq:laxstrict}) to obtain the (unique, up to a 
space translation) profiles $z_{0}$ and $z_{n}$ corresponding to the intervals
$[z_{+}, z^{*}_{1}]$ and $[z^{*}_{2n-1},z_{-}]$. 
Taking into account also the results of Lemma \ref{lem:intersection}, 
we know that such maximal solutions verify, up to a space 
translation, the following properties:
\begin{itemize}
    \item[(i)] $z_{0}: [0, +\infty) \to (z_{+}, z_{0}(0)]$, $z_{0}$ 
    monotone decreasing, $z'_{0}$ monotone increasing and 
    $z_{0}(+\infty) = z_{+}$;

    \item[(ii)] $z_{n}: (-\infty, 0] \to [z_{n}(0), z_{-})$, 
    $z_{n}$  and $z'_{2n}$
    monotone decreasing and $z_{n}(-\infty) = z_{-}$;

    \item[(iii)] $z_{0}(0)\geq z^{*}_{1}$ and $z_{n}(0)\leq z^{*}_{2n-1}$.
\end{itemize}
It is worth to observe that the existence of the above profiles with 
the aforementioned properties does not depend on the value of $\e>0$.
Hence, we are left  with the construction of the \emph{intermediate maximal solutions} 
$z_{k,l}$ and $z_{k,r}$,  $k=1,\ldots,n-1$, solutions 
respectively of 
\begin{equation}
    \e z'' = z-h_{2k}(z')
    \label{eq:crescente}
\end{equation}
and
\begin{equation}
    \e z'' = z-h_{2k+1}(z'),
    \label{eq:decrescente}
\end{equation}
where $h_{i}$ denotes the inverse of $F(\cdot;s)$ on $I_{i}$,  $i=1,\ldots,2n$.

\begin{proposition}\label{prop:intermediate}
    Let us assume conditions (\ref{eq:ole2}) and (\ref{eq:laxstrict}) hold. Then, for any 
    $k=1,\ldots,n-1$, there exists a (unique up to space translations) maximal solution 
    $z_{k,l}$ of (\ref{eq:crescente}) and $z_{k,r}$ of (\ref{eq:decrescente}), with initial data
    \begin{equation}
	\begin{cases}
	    z_{k,l}(0) =z_{k,r}(0) 
		=z_{2k}^{*}&   \\
		z'_{k,l}(0) =z'_{k,r}(0) 
		    =F(z_{2k}^{*};s). & 
	    \label{eq:data}
	\end{cases}
    \end{equation}
    Moreover
    \begin{equation}
	z''_{k,l}(0) = z''_{k,r}(0) = 0
	\label{eq:continua}
    \end{equation}
    and
    \begin{itemize}
	\item[(i)] $z_{k,l}$, $z'_{k,l}$ are 
    monotone decreasing and $z_{k,l}$ is not globally defined, that 
    is there exists a point $\bar\xi_{k,l}>0$ such that
    \begin{equation*}
	z_{k,l}(\bar\xi_{k,l}) - \e z''_{k,l}(\bar\xi_{k,l}) = z^{*}_{2k-1},\quad 
	z'_{k,l}(\bar\xi_{k,l}) = F(z^{*}_{2k-1};s);
    \end{equation*}
    
	\item[(ii)] $z_{k,r}$ is monotone decreasing,  $z'_{k,r}$ 
    monotone increasing and $z_{k,r}$ is not globally defined, that 
    is there exists a point $\bar\xi_{k,r}<0$ such that
    \begin{equation*}
	z_{k,r}(\bar\xi_{k,r}) - \e z''_{k,r}(\bar\xi_{k,r}) = z^{*}_{2k+1},\quad 
	z'_{k,r}(\bar\xi_{k,r}) = F(z^{*}_{2k+1};s);
    \end{equation*}

	\item[(iii)] there exists an $\e_{k}$ such that, for any 
	$\e<\e_{k}$, 
	\begin{equation}
	    z_{k,l}(\bar\xi_{k,l})\leq z^{*}_{2k-1},\quad 
	    z_{k,r}(\bar\xi_{k,r})\geq z^{*}_{2k+1}.
	    \label{eq:intersection2}
	\end{equation}
    \end{itemize}
\end{proposition}

\begin{proof}
    Let us start by justifying the choice of initial data in 
    (\ref{eq:data}). Since we want to joint $z_{k,l}$ and $z_{k,r}$ 
    in $0$ to obtain a smooth (say, $C^{2}$) profile, we are forced to choose 
    (\ref{eq:data})$_{2}$. Moreover, differentiating  
    (\ref{eq:original}) we get
    \begin{equation}
	z'' = (z-\e z'')' F'(z-\e z'';s).
	\label{eq:second}
    \end{equation}
    In the sequel, we shall also need $z'''_{k,l}(0)$ and $z_{k,r}'''(0)$ 
    finite and therefore, assuming that requirement, form (\ref{eq:second}) we obtain 
    (\ref{eq:continua}), which implies (\ref{eq:data})$_{1}$, in view 
    of (\ref{eq:crescente}) and (\ref{eq:decrescente}).
    
    Since $F(z^{*}_{2k};s)<0$, we can extend 
	continuously the function $h_{2k}$ at the right of 
	$z^{*}_{2k}$ to obtain the existence of a (not necessary 
	unique) $C^{2}$ maximal solution of 
    (\ref{eq:crescente})--(\ref{eq:data}), which coincides with 
    $z_{k,l}(\xi)$, for any $\xi>0$,  as long as $z_{k,l}(\xi)\leq 
    z_{2k}^{*}$. Hence, property (i) is proved as before: $h_{2k}$ is 
    the inverse of $F(\cdot;s)$ on an interval where $F(\cdot;s)$ is 
    increasing and therefore it corresponds to the case studied in 
    Proposition \ref{prop:maximalplus}.
    A symmetric argument will lead to the situation of  Proposition \ref{prop:maximalminus}, 
    which gives (ii).
    The remaining part of the results will be proved only for 
   $z_{k,l}$, the ones for $z_{k,r}$ being similar and thus we drop 
   the subscript $k,l$.
   
    Let us now prove
    \begin{equation}
	|z'''(0)|<+\infty.
	\label{eq:bound}
    \end{equation}
    To this end, let us come back to the original variables $u$ and $q$. From 
    (\ref{eq:nonconvex}), and after integration with respect to $\xi$ of 
    the first line, the equations for the profiles $u$ and $q$ are given by
    \begin{equation*}
	\begin{cases}
	    F(u;s) =- q  & \\
	    u'=\e q''-q,
	\end{cases}
    \end{equation*}
    because $q(\pm\infty) =0$. Thus,  the  profile $u$, regular away from $0$ and 
    continuous there, verifies
    \begin{equation*}
	u' =- \e F''(u;s) + F(u;s) = - \e F'(u;s)u''
	- \e F''(u;s)(u')^{2} + F(u;s).
    \end{equation*}
    Hence, we obtain the following system in the state space $(u,v=u')$
    \begin{equation}
	\begin{cases}
	    u'= v & \\
	    \e F'(u;s)v'=-\e F''(u;s)v^{2}-v+F(u;s).
	\end{cases}
	\label{eq:uvprofile2}
    \end{equation}
    At this point, as in \cite{KN99}, we shall introduce a new 
    independent variable $\eta$ given by
    \begin{equation*}
	\xi = \int_{\eta}^{\infty}F'(u(\zeta);s)d\zeta,
    \end{equation*}
    which will move the singularity attained for $\xi =0$ to $\eta= 
    \infty$, thanks to exponential decay of $F'(u(\eta);s)$ toward zero, 
    as $\eta\to \infty$ (the same kind of procedure used for 
    regularity results in Section \ref{sec:reg}). Hence, 
    (\ref{eq:uvprofile2}) becomes
    \begin{equation}
	  \begin{cases}
	      u'= F'(u;s)v & \\
	      v'=\frac{1}{\e}\left (-\e F''(u;s)v^{2}-v+F(u;s)\right ).
	  \end{cases}
	  \label{eq:uvprofile3}
      \end{equation}
      The singular point $u = z^{*}_{2k}$ for (\ref{eq:uvprofile2}) 
      corresponds to a pair of equilibrium points $(z^{*}_{2k}, v^{*}_{1,2})$ 
      for (\ref{eq:uvprofile3}),
      with
      \begin{equation*} 
	  v^{*}_{1,2} = \frac{-1 \pm \sqrt{1+
	 4\e |F''(z^{*}_{2k};s)F(z^{*}_{2k};s)|}}{2\e F''(z^{*}_{2k};s)}.
	\end{equation*}
    Moreover, from (i) we know that
    \begin{equation*}
	v=u' = z'-\e z''' = \frac{z''}{\e F'(z-\e z'';s)}<0.
    \end{equation*}
    Hence the orbit  lies in the lower 
    halfplane $v<0$. Moreover,   $u\to z_{2k}^{*}$ as $\eta\to 
    -\infty$, in view of the requirements (\ref{eq:data}) 
    in the original variables and a straightforward analysis of 
    the dynamical system   (\ref{eq:uvprofile3}) 
    ensures 
    the   $(z^{*}_{2k}, v^{*}_{1})$ is a saddle point (see also Section \ref{sec:reg} for further 
    details on a very close dynamical system).   
    Thus there exists an unique orbit for which $u(-\infty) = z^{*}_{2k}$, 
    namely the unique orbit ending at the saddle point. 
    In particular, the value $v$ converges toward the finite value $v^{*}_{1}<0$, as 
    $\eta \to - \infty$. Hence, coming back to 
    the original independent variable $\xi$,  $u'$ is bounded in $0$, 
    that is (\ref{eq:bound}). Moreover, as already pointed out 
    before, that property implies $z''(0) =0$. 
    
    We are now left with the proof of (iii), namely $z(\bar\xi)\leq z^{*}_{2k-1}$. 
    As for Lemma \ref{lem:intersection}, we shall prove that property 
    by comparing in the phase space $(z,z')$, the graph $\varphi(z)$ 
    of the trajectory of $z(\xi)$ with the one $\psi(z) = F(z;s)$ of $z' = F(z;s)$ close 
    to the point $(z^{*}_{2k},F(z^{*}_{2k};s))$.
   Since $ \psi'(z^{*}_{2k}) = F'(z^{*}_{2k};s)=0$ and
   \begin{equation*}
       \varphi'(z^{*}_{2k}) = \left . \frac{z''}{z'} \right |_{\xi=0} 
       =0, 
   \end{equation*}
     we need to 
   evaluate the second derivative at this point. Clearly, 
   $\psi''(z^{*}_{2k}) = F''(z^{*}_{2k};s)$, while
   \begin{equation*}
       \varphi''(z^{*}_{2k}) = \left. \frac{\displaystyle{\frac{d}{d\xi} \left ( 
		\frac{z''}{z'}  \right )} }{z'}\right | _{\xi =0} = 
       \frac{z'''(0)z'(0) - z''(0)^{2}}{(z'(0))^{3}} = 
       \frac{z'''(0)}{F(z^{*}_{2k};s)^{2}}.
   \end{equation*}
   Moreover, from the above calculations, we know that
   \begin{equation*}
       z'''(0) = \frac{z' (0) - u'(0)}{\e} = \frac{F(z^{*}_{2k};s) - 
       v^{*}_{1}}{\e}.   
   \end{equation*}
   Therefore, $z'''(0)$ is the biggest solution of the equation
   \begin{equation*}
       Q(\lambda) := (\e \lambda - F(z^{*}_{2k};s))^{2} - 
       \frac{\lambda}{F''(z^{*}_{2k};s)} = 0.
   \end{equation*}
   Since
   \begin{equation*}
       Q(F(z^{*}_{2k};s)^{2}F''(z^{*}_{2k};s)) = \e
       F(z^{*}_{2k};s)^{3}F''(z^{*}_{2k};s)(\e F(z^{*}_{2k};s)F''(z^{*}_{2k};s) -2),
   \end{equation*}
   and $F(z^{*}_{2k};s)F''(z^{*}_{2k};s)>0$, there exists an 
   $\e_{0}>0$ such that, for any $\e <\e_{0}$, 
   $Q(F(z^{*}_{2k};s)^{2}F''(z^{*}_{2k};s)) <0$, that is
   \begin{equation*}
       F(z^{*}_{2k};s)^{2}F''(z^{*}_{2k};s) < z'''(0),
   \end{equation*}
   which implies 
   \begin{equation*}
      0> \varphi''(z^{*}_{2k}) > \psi''(z^{*}_{2k}).
   \end{equation*}
   Therefore, the graph $\varphi(z)$ leaves the point $z^{*}_{2k}$ 
   above $\psi(z)$, while going to the left. Hence, we can argue as in 
   Lemma \ref{lem:intersection} to conclude $z(\bar\xi) \leq 
   z^{*}_{2k-1}$ and the proof is complete.
\end{proof}

In view of the above result, as in Theorem \ref{theo:existence}, we can glue together 
the profiles $z_{0}$, $z_{k,l}$, $z_{k,r}$, $k=1,\ldots,n-1$, $z_{n}$ 
to obtain the desired radiating profile, joining $u_{-}=z_{-}$ and $u_{+}=z_{+}$.

\begin{proposition}\label{theo:existence2}
Under conditions (\ref{eq:ole2}) and (\ref{eq:laxstrict}), there exists an $\bar \e>0$ 
such that, for any $\e<\bar\e$,  there exists a (unique  up to space translations)  
$C^{1}$ profile $z$ with $z(\pm\infty) = z_{\pm}$ and a speed $s$ such that the function 
$z(x-st)$ is solution of (\ref{eq:original}). 
This solution is  $C^{2}$ away from the $n$ points $z^{*}_{2k-1}$, $k=1,\ldots,n$, where $z''$ 
has at most a jump discontinuity.
   
   Moreover, there exists a (unique up to space translations)  
   profile $u$  with $u(\pm\infty) = u_{\pm}$ and a speed $s$ (given by the Rankine--Hugoniot
   condition) such that the function $u(x-st)$ is solution of (\ref{eq:nonconvexprof}). 
   This profile is continuous away from the $n$ points $z^{*}_{2k-1}$, $k=1,\ldots,n$,
   where it has at most a jump discontinuity.
   At these points, the Rankine--Hugoniot and the admissibility conditions of the scalar 
   conservation law $u_{t}+f(u)_{x}=0$ are satisfied.
\end{proposition}

\begin{proof}
   First of all, from Proposition \ref{prop:intermediate},  we can 
   define $C^{2}$ profiles $z_{1}, \ldots, z_{n-1}$, gluying together 
   the profiles $z_{k,l}$, $z_{k,r}$, after an 
   appropriate space translation. Hence, we end up with $n+1$ profiles 
   $z_{0},\ldots,z_{n}$, all of them decreasing.
   
   Moreover, for $\e<\bar\e = \min\{\e_{k}\}_{k=1,\ldots,n-1}$, there 
   exist points in the 
    $(z,z')$ plane where the graphs of two consecutive $z_{k}$ can be 
    glued together in that point, in order to give a $C^{1}$ solution 
    of (\ref{eq:original}).
    Moreover, due to the monotonicity of these graphs, again 
    consequence of Proposition \ref{prop:intermediate}, this intersection is indeed unique.  
    Hence 
    with  appropriate space translations, we can find a points $\bar 
    \xi_{k}$, $k=1,\ldots,n-1$,
    such that $z_{k}(\bar \xi_{k}) = z_{k+1}(\bar \xi_{k}) = \bar 
    z_{k}$ and $z'_{k}(\bar \xi_{k}) = 
    z'_{k+1}(\bar \xi_{k})= \tilde z_{k}$, $k=1,\ldots,n-1$.
    Now we define 
    \begin{equation*}
	z(\xi) = 
	\begin{cases}
	    z_{n}(\xi) &  \xi\in (-\infty , \bar \xi_{n-1}] \\
	     z_{k}(\xi),    & \xi\in[\bar\xi_{k},\bar\xi_{k-1}], \qquad
	     k=2,\ldots,n-1, \\
	    z_{0}(\xi) & \xi\in[\bar\xi_{0},+\infty).
	\end{cases}
	\end{equation*}
    Thus, this profile defines a $C^{1}$ solution of
    (\ref{eq:original}) and it verifies $z(\pm\infty) = z_{\pm}$. 
    Moreover, $z(x)$ is $C^{2}$
    away from the points $\bar\xi_{k}$, where it verifies
    $\e z''(\bar \xi_{k} -0) = \e z''_{k+1}(\bar \xi_{k}) 
    =\bar z_{k} - h_{2k+2}(\tilde z_{k})$ and
    $\e z''(\bar \xi_{k} +0) = \e z''_{k}(\bar \xi_{k}) = 
    \bar z_{k} - h_{2k+1}(\tilde z_{k})$.
   
   The regularity of $u=u(x-st)$ is a direct consequence of 
   the first part of the theorem and of the relation $u=z-\e z''$. 
   Moreover, in the case of a discontinuity in $u$, namely $u(\bar \xi_{k} 
   -0)\neq u(\bar \xi_{k} +0)$,  $(u(\bar \xi_{k} -0), u(\bar \xi_{k} +0); s)$
   verifies 
      the Rankine--Hugoniot condition 
      for the strictly 
	 convex conservation law $u_{t}+f(u)_{x}=0$. Indeed, $u(\bar \xi_{k} -0) = 
   h_{2k+2}(\tilde z_{k})$, $u(\bar x +0) = h_{2k+1}(\tilde z_{k})$ and a direct 
   calculation shows
   \begin{equation*}
       \frac{f(h_{2k+1}(\tilde z_{k}))-f(h_{2k+2}(\tilde z_{k}))}{ 
       h_{2k+1}(\tilde z_{k})- h_{2k+2}(\tilde z_{k})} = 
	s.
   \end{equation*}
   In addition, $u(\bar \xi_{k} -0)  = h_{2k+2}(\tilde z_{k})\in 
   I_{2k+2}$,  $u(\bar \xi_{k} +0)  = h_{2k+1}(\tilde z_{k})\in 
   I_{2k+1}$ and $F(u(\bar \xi_{k} -0);s) = F(u(\bar \xi_{k} +0);s) = \tilde z_{k}$. 
   Since $F(\cdot;s)$ has a local minimum and no local 
   maximum in $(u(\bar \xi_{k} +0),u(\bar \xi_{k} -0))$, we conclude 
   $F(u;s) < F(u(\bar \xi_{k} \pm 0);s)$, for any $u\in (u(\bar 
   \xi_{k} +0),u(\bar \xi_{k} -0))$, that is the Oleinik 
   condition for the inviscid shock 
   $(u(\bar \xi_{k} -0), u(\bar \xi_{k} +0); s)$.
\end{proof}

Thanks to the above theorem, we know there exist radiative shocks, if 
Oleinik coindition (\ref{eq:ole2}) and the non degeneracy 
(\ref{eq:laxstrict}) are satisfied. Hence, we can pass now 
 to the proof of qualitative  properties of that solutions. 
 In particular,  for a 
profile $u$ that belongs to $BV$, we shall prove monotonicity and uniqueness in the class 
modulo $L^{1}$. It is worth to observe that these properties do not require 
assumption (\ref{eq:laxstrict}), because they are based on the contraction 
properties of the model under consideration. Actually,  they 
 will be used in the existence result of Theorem 
 \ref{theo:generalflux} in the case of general smooth flux functions 
 which may violate (\ref{eq:laxstrict}).
 
\begin{theorem}\label{theo:uproperties}
    Let $(u,q)$ be a radiative shock solution of (\ref{eq:nonconvex}). 
    If $u\in BV$, then $u$ is monotone. 
    Moreover, given two $BV$ radiative shocks $(u,q_{1})$ and 
    $(v,q_{2})$, such 
    that $u-v \in L^{1}$, then $u$ and $v$ are equal up to a space translation.
\end{theorem}

\begin{proof} 
As pointed out before by several authors, \cite{KNN98, 
    KN99b,KN99,Ser03,Ser04,LM03}, it is convenient to 
    rewrite the $2\times 2$ system (\ref{eq:nonconvex}) as follows:
    \begin{equation}
        u_{t} + f(u)_{x} = -\frac{1}{\e} \left (u - K^{\e}* u \right 
	),
        \label{eq:rewritten}
    \end{equation}
    where the convolution kernel $K^{\e}$ is given by
    \begin{equation*}
        K^{\e}(x) = \frac{1}{2\sqrt\e}e^{-\frac{|x|}{\sqrt\e}}.
    \end{equation*}
Let $u$ and $v$ be radiative profiles of class $BV$, associated  
 with a shock triplet $(u_-,u_+;s)$. Assume that $v-u$ is  
 integrable, a fact that is certainly true if $u$ is in $BV$ and $v$  
 is a shift $\tau_hu$ of $u$. Finally, let us denote $w:=v-u$. We  
 integrate the entropy inequality
 $$
  |w|_t+\Bigl((f(v)-f(u)){\rm sgn}(w)\Bigr)_x\le
  \frac 1\e\,(K^\e*w-w){\rm sgn}(w)
 $$
 and obtain, with an argument {\it \`a la Kruzhkov},
 $$
 \e\,\frac{d}{dt} \,\|w\|_1\le\int_{\R} (K^\e w){\rm sgn}(w) \,dx-\|w\|_1.
 $$
 Since $w$ is a function of $x-st$, the norm $\|w\|_1$ is constant  
 and the left-hand side above is zero. Therefore, we have
 $$
 \|w\|_1\le\int_{\R} (K^\e*w){\rm sgn}(w)\,dx.
 $$
 Since $K^\e$ has positive values and unit integral, we deduce that $w$  
 has a constant sign.

 Let $h$ be a real number. Since $\tau_hv:=v(\cdot+h)$ is also a  
 profile and $\tau_hv-v\in L^1(\R)$, the above fact may be applied  
 to the pair $(\tau_hv,u)$. We obtain that for every real number $h 
 $, $\tau_hv-u$ has a constant sign.

 Let us focus on the case $v=u$. Since the integral of $\tau_hu-u$  
 equals $h(u_+-u_-)$, we see that the sign of $\tau_hu-u$ is that of  
 $h(u_+-u_-)$. In other words, $u$ is monotonous.

 Going back to the general case, the integral of $\tau_hv-u$ equals  
 $h(u_+-u_-)$ plus a constant (the integral of $v-u$). Hence there  
 exists a number $h$ for which this integral equals zero. But since $ 
 \tau_hv-u$ has a constant sign, this means that $\tau_hv-u\equiv0$  
 almost everywhere.
\end{proof}

Now we are ready to remove the non degenerate assumption 
(\ref{eq:laxstrict}) in the existence of a radiative shock. 
We perform this task by means of an 
approximationg procedure already used in \cite{Ser06} for discrete 
shock profiles for conservation laws.
\begin{theorem}\label{theo:generalflux}
    Under conditions (\ref{eq:ole2}),
       there exists an $\bar \e>0$ such that, for any $\e<\bar\e$, 
      there exists a (unique  up to space translations)   
      $C^{1}$ profile  $z$ with $z(\pm\infty) = z_{\pm}$ such that the function $z(x-st)$ 
      is a solution of (\ref{eq:original}),
where the speed $s$ is given by the Rankine--Hugoniot condition. 
This solution is  $C^{2}$ away from the 
      $n$ points $z^{*}_{2k-1}$, $k=1,\ldots,n$, where $z''$ has at most a jump discontinuity.

Moreover, there exist a  (unique  up to space translations)
      profile $u$ with $u(\pm\infty) = u_{\pm}$ and
a speed $s$ such that the function $u(x-st)$ is solution of (\ref{eq:nonconvexprof}),
where the speed $s$ is given by the Rankine--Hugoniot condition. 
This profile is continuous away from the $n$ points $z^{*}_{2k-1}$, $k=1,\ldots,n$,
 where it has at most a jump discontinuity which verifies the 
Rankine--Hugoniot and the admissibility conditions of the scalar 
conservation law $u_{t}+f(u)_{x}=0$.
\end{theorem}

\begin{proof}
      We start by approximating  a $C^{2}$ flux function $f$ satisfying 
      (\ref{eq:ole2}) with a sequence of smooth functions $f^{n}$ 
      satisfying both (\ref{eq:ole2}) and (\ref{eq:laxstrict}). Then, 
      given an inviscid shock $(u_{-},u_{+};s)$, $u_{-}>u_{+}$, 
      Proposition \ref{theo:existence2} gives the existence of radiative 
      shocks $(u^{n},z^{n})$ for any $n$ and for any $\e< \bar\e$ 
      such that $u^{n}(\pm\infty) = z^{n}(\pm\infty) = u_{\pm}$. 
      This  value $\bar\e $ is independent 
      from $n$, because   $f^{n}$ and 
      $(f^{n})''$ remain bounded,  (see the constraints that gives the values $\e_{k}$ in 
      Proposition \ref{prop:intermediate} and hence $\bar\e$ in 
      Proposition \ref{theo:existence2}). We fix these profiles with the 
      condition
      \begin{equation*}
          \frac{1}{2}(u_{-}+u_{+})\in [u^{n}(0+),u^{n}(0-)].
      \end{equation*}
      Since $u^{n}$ is decreasing,
      \begin{equation*}
          \|\tau_{h} u^{n} -u^{n}\|_{L^{1}} = |h(u_{-}-u_{+})|,
      \end{equation*}
      which implies that $u^{n}$ is equicontinuous 
      in $L^{1}$.
      Moreover, $u_{+}\leq u^{n} \leq u_{-}$, that is, $u^{n}$ is 
      equibounded in $L^{\infty}$ and hence in $L^{1}_{loc}$. We 
      recall  that the solution $z^{n}$ of $-\e 
      (z^{n})'' + z^{n} = u^{n}$ with $z^{n}(\pm\infty)=u_{\pm}$ is given by
      \begin{equation*}
          z^{n}(\xi) = \frac{1}{2\sqrt\e}\int_{-\infty}^{+\infty} 
	  e^{-\frac{|\zeta|}{\sqrt\e}} u^{n}(\xi-\zeta)d\zeta
      \end{equation*}
      and it satisfies the same estimates. Thus, passing if necessary 
      to subsequences, we have $u^{n}\to u$ and $z^{n}\to z$ in 
      $L^{1}_{loc}$ and bounded almost everywhere. Clearly, $u$ is 
      decreasing and therefore 
      \begin{equation*}
          u_{+} \leq u(+\infty) \leq \frac{1}{2}(u_{-}+u_{+})  \leq u(-\infty) \leq u_{-}.
      \end{equation*}
      Then, we pass to the limit in the profile equations
      \begin{equation}
	  \begin{cases}
	      F^{n}(u^{n};s) =(z^{n})' & \\
	      -\e (z^{n})'' + z^{n} =u^{n}, &
	  \end{cases}
	  \label{eq:nonconvexprofn}
      \end{equation}
      where $F^{n}(u;s) = f^{n}(u)-f^{n}(u_{\pm}) - s(u-u_{\pm})$ 
      to conclude that $(u,z)$ defines a radiative profile for the 
      inviscid shock $(u(-\infty),u(+\infty);s)$. Hence we are left 
      to the proof of $u(\pm\infty) = u_{\pm}$. 
      Integrating (\ref{eq:nonconvexprofn})$_{1}$ and passing into 
      the limit a.e.\ we get
      \begin{equation*}
          z(\xi) - z(\zeta) = \int_{\zeta}^{\xi} F(u(x);s)dx.
      \end{equation*}
      Since $z(\pm\infty) = u(\pm\infty)\in\R$ and $u\in 
      [u_{+},u_{-}]$, from the above 
      relation and from (\ref{eq:ole2}) we conclude
      \begin{equation*}
         0 >  \int_{-\infty}^{+\infty} F(u(x);s)dx > -\infty.
      \end{equation*}
      Therefore in particular $F(u(\pm\infty);s)=0$ and, using 
      (\ref{eq:ole2}) and
      \begin{equation*}
          u(+\infty) \in \left [u_{+}, \frac{1}{2}(u_{-}+u_{+}) 
	  \right ], \ u(-\infty) \in \left [, 
	  \frac{1}{2}(u_{-}+u_{+}) u_{-} \right ],
      \end{equation*}
      we conclude $u(\pm\infty) = u_{\pm}$.
\end{proof}

\begin{remark}\label{rem:regularitynonconvex}
    Proceeding as in Section \ref{sec:reg}, it is possible to prove 
    the radiative shocks of (\ref{eq:nonconvex})  increase their regularity 
    as $\e\downarrow 0$ also in 
    this general non convex case. Indeed, across the values $z^{*}_{2k-1}$, 
    $k=1,\ldots,n$, where   $F(\cdot;s)$ has local minima, we are in 
    the same 
    situation of  the one discussed in the convex case 
    and the regularity increases as $\e\downarrow 0$. On the contrary, across the values $z^{*}_{2k}$, 
    $k=1,\ldots,n$, where   $F(\cdot;s)$ has local maxima, we can 
    repeat the arguments of Proposition \ref{prop:reg} and we obtain 
    that the 
    trajectory in the saddle point corresponding to the value 
    $u=z^{*}_{2k}$  is never tangent to the eigenvector $(0,1)^t$
   of the negative eigenvalue and therefore the profile is regular. 
\end{remark}

\section{Reduction from the system case to the scalar case}\label{sec:system}

Let us consider the following strictly hyperbolic--elliptic coupled system
\begin{equation}
    \begin{cases}
	u_{t}+ f(u)_{x} +Lq_{x}=0 & \\
	-q_{xx} + Rq +G\cdot u_{x} =0, &
    \end{cases}
    \label{eq:system}
\end{equation}
where $x\in\R$, $t>0$, $u\in\R^{n}$, $q$ is scalar and $R>0$, 
$G$, $L\in\R^{n}$ are constants. 
In this section we shall prove that the existence of travelling wave
solutions of (\ref{eq:system}), for sufficiently small shocks,  
reduces to the study of a scalar
model,  and therefore we shall obtain their
existence (and regularity) as corollary of the previous sections.

According to Definition \ref{def:radiative}, let us consider a radiative shock   wave solutions for
(\ref{eq:system}), namely a solution of the form
$(u,q) = (u(x-st),q(x-st))$, with $s$ verifying the Rankine--Hugoniot condition
\begin{equation*}
    u(\pm\infty) = u_{\pm},\quad s(u_{+}-u_{-}) = f(u_{+})-f(u_{-}).
\end{equation*}
As usual, we denote with $\lambda_{1}(u) < \dots <\lambda_{n}(u)$ the 
$n$ real eigenvalues of the matrix $\nabla f(u)$. 
Neglecting the higher order term $-q_{xx}$ in (\ref{eq:system})$_{2}$ 
(see, for instance \cite{KNN98,LM03}), we end up with the reduced
system
\begin{equation*}
    u_{t} + f(u)_{x} = R^{-1}G\cdot u_{xx}L =R^{-1} L\otimes G u_{xx}
\end{equation*}
and the $k$th field has a diffusive dynamics  near $u_{\pm}$, provided
\begin{equation}
    \ell_{k}(u_{\pm})\cdot L\otimes G r_{k}(u_{\pm}) = G\cdot r_{k}(u_{\pm})
    \ell_{k}(u_{\pm})\cdot L >0,
    \label{eq:mainassump}
\end{equation}
where, as usual,  $r_{k}(u)$ and $\ell_{k}(u)$ denote the $k$th right and
left eigenvector of $\nabla f(u)$, normalized such that $\ell_{i}(u)\cdot
r_{j}(u)=\delta_{ij}$.

Once again, we   introduce the variable $z$ as the opposite of the 
antiderivative of $q$, that
is $-z_{x}:=q$ and therefore $z(\pm\infty)=z_{\pm}=G\cdot u_{\pm}$.
Hence, we proceed as in the previous sections to conclude the  
the system for  $u=u(x-st)$ and $z=z(x-st)$ is given by
\begin{equation}
    \begin{cases}
	Lz'=   f(u)-f(u_{\pm}) -s(u-u_{\pm})  & \\
      Rz -z'' =G\cdot u. &
    \end{cases}
    \label{eq:systemprof}
\end{equation}
Given the vector $L$, let $P:\R^{n}\to\R^{n}$ and $Q:\R^{n}\to\R$ be 
linear applications such that
\begin{equation*}
     \ker P = \textrm{span}\,\{L\},\qquad QL =1.
\end{equation*}
Then system (\ref{eq:systemprof}) becomes
\begin{equation}
    \begin{cases}
	P\left (f(u)-f(u_{\pm}) -s(u-u_{\pm}) \right )= 0 & \\
	z' = Q\left (f(u)-f(u_{\pm}) -s(u-u_{\pm}) \right ) & \\
	G\cdot u =Rz- z''. &
    \end{cases}
    \label{eq:systemprof2}
\end{equation}
Hence, the existence of our profile reduces to a scalar case,
provided the   constraints (\ref{eq:systemprof2})$_{1}$, together 
with condition (\ref{eq:mainassump}), ensure
\begin{equation}\label{eq:reduction}
    Q\left (f(u)-f(u_{\pm}) -s(u-u_{\pm}) \right ) = \hat F(G\cdot u;s),
\end{equation}
  in a neighborhood of $(u_{\pm};\lambda_{k}(u_{\pm}))$, that 
is for $|u_{+}-u_{-}|$ sufficiently small.
This property is proved in  the next lemma.

\begin{lemma}\label{lem:reduction}
 Let $f$ be $C^{1}$ and assume condition (\ref{eq:mainassump}) holds. 
Then, there exists  a neighborhood $\mathcal U \times I$ of $(u_{\pm};\lambda_{k}(u_{\pm}))$
such that for any $(u;s)\in{\mathcal U\times I}$,
\begin{equation*}
	P\left (f(u)-f(u_{\pm}) -s(u-u_{\pm}) \right )=0 \ \iff\ u=\Phi(G\cdot u;s), 
\end{equation*}
with $\Phi\in C^1(\R\times \R;{\mathcal U})$, not depending on the choice
of $P$ such that $\ker P = \textrm{span}\,\{L\}$.
\end{lemma}

\begin{proof}
     Assume, without loss of generality, 
     $(u_{\pm};\lambda_{k}(u_{\pm}))=0$ and $\|G\|=1$.
    Let $V_1,\dots, V_{n-1}$ be a basis for $G^\perp$, so that
     $\{V_1,\dots, V_{n-1},G\}$ is a basis for $\R^n$.
    For $u= c_1V_1+\cdots+c_{n-1}V_{n-1}+\alpha G$, let us denote with
     $
       \phi(c_1,\dots,c_{n-1},\alpha;s)
     $
     the function $P\left (f(u)-f(u_{\pm}) -s(u-u_{\pm}) \right )$. 
    Then $\phi(0,\dots,0;0)=0$.
    Moreover, for $j\in\{1,\dots,n-1\}$,
    \begin{equation*}
       \frac{\partial\phi}{\partial c_j}(0) =P\,\nabla f(0)\,V_j.
       \end{equation*}
    Then the conclusion follows from the Implicit Function Theorem,
    provided
    the vectors $P\nabla f(0)V_1,\dots,P\nabla f(0)V_{n-1}$,
    $j=1,\ldots,n-1$ are linearly independent, because, in that case, 
    the relation
    $\phi(c_1,\dots,c_{j-1},\alpha; s)=0$ can be locally written as a 
    function of $(\alpha=G\cdot U ; s)$.

    Let $a_1,\dots, a_{n-1}\in\R$ be such that 
    \begin{equation*}
	a_1P\nabla f(0)V_1+\dots+a_{n-1}P\nabla f(0)V_{n-1}=0.
	\end{equation*}
    Then 
    \begin{equation*}
       P\nabla f(0)\left(a_1V_1+\dots+a_{n-1}V_{n-1}\right)=0
   \end{equation*}
    and, by definition of $P$, there exist $\beta\in\R$ such that
    \begin{equation*}
	\nabla 
	f(0)\left(a_1\,V_1+\dots+\alpha_{n-1}\,V_{n-1}\right)=\beta L.
    \end{equation*}
    Applying $\ell_{k}(0)$ at the left hand side, we have
    \begin{equation*}
       0=\ell_{k}(0) \nabla 
       f(0)\left(a_1V_1+\dots+a_{n-1}V_{n-1}\right)=\beta \ell_{k}(0)\cdot L.
   \end{equation*} 
    Hence, condition (\ref{eq:mainassump}) implies $\beta=0$, namely
      $a_1V_1+\dots+a_{n-1}V_{n-1}\in~\ker~\nabla f(0)$,  which means  
    \begin{equation*}
	a_1V_1+\dots+a_{n-1}V_{n-1} = \gamma 
	    r_{k}(0),
    \end{equation*} 
    for some $\gamma\in\R$. Since  $V_1,\dots, V_{n-1}$ is 
    a basis for $G^\perp$, using again (\ref{eq:mainassump}) we 
    conclude $\gamma =0$ and therefore $a_1=\dots= a_{n-1}=0$.
    This means that  $P\nabla f(0)V_j$, $j=1,\ldots,n-1$ are linearly independent
    and Implicit Function Theorem cann be applied.
    
 Finally, let us note that since $\ker P = \textrm{span}\,\{L\}$,
 $P\left (f(u)-f(u_{\pm}) -s(u-u_{\pm}) \right )=0$ is equivalent to
$f(u)-f(u_{\pm}) -s(u-u_{\pm})\in\textrm{span}\,\{L\}$, hence the function $\Phi$
does not depend on the specific choice of $P$.
  
The proof is complete. 
\end{proof}

At this point, we assume the $k$th field is genuinely nonlinear, that is
\begin{equation}\label{eq:GNL}
    \nabla \lambda_{k}(u) \cdot r_{k}(u)   \neq 0
\end{equation}
and let $(u_{-},u_{+};s)$ be a $k$--Lax radiating shock for (\ref{eq:system}), that is
\begin{equation*}
    \lambda_{k}(u_{+}) < s < \lambda_{k}(u_{-}), \quad 
    \lambda_{k-1}(u_{-}) < s < \lambda_{k+1}(u_{+}).
\end{equation*}
We shall prove that in this framework, in the reduced scalar dynamics yielding
the existence of our profile, that is (\ref{eq:systemprof2})$_{2}$, the flux function $F(\cdot;s)$  
does not change convexity. 
On the other hand, in  the general case and assuming the states verify the Liu E--condition, 
the reduction will lead to the non convex model treated in Section \ref{subsec:nonconvex}.

\begin{proposition}\label{prop:convexreduction}
Let us assume condition (\ref{eq:mainassump}) and (\ref{eq:GNL}) hold. 
Then the function $\hat F(\cdot;s):\R\to\R$ defined in  (\ref{eq:reduction}) 
is either strictly convex or strictly concave in a neighborhood of $z_{\pm}=G\cdot u_{\pm}$. 
\end{proposition}

\begin{proof}
   Since $(u_{-},u_{+};s)$ form a $k$--Lax shock for
(\ref{eq:system}), it suffices to prove 
\begin{equation*}
   \frac{d^{2}}{dw^{2}}\hat F(G\cdot u_\pm; s)\neq 0,
\end{equation*}
for sufficiently small shocks.

Since $F(u;s):= f(u) -  f(u_{\pm}) - s(u-u_{\pm})=0$,
there hold $P\,F(u;s) =0$ and $Q\,F(u;s) = \hat F(G\cdot u;s)$ for any $u\in\mathcal{U}$.
Hence
	  \begin{equation*}
	        F(u;s) =  \hat F(G\cdot u;s)\,L.
	  \end{equation*}
Differentiating with respect to $u$, one has
	  \begin{equation}
		\nabla   f(u) - s I   =  \frac{d}{dw} \hat F(G\cdot u;s)\,L\otimes G
	      \label{eq:differ}
	  \end{equation}
Applying $\ell_{k}(u)$ and $r_{k}(u)$ respectively to the left and to the right of (\ref{eq:differ}), we get
\begin{equation}
  \lambda_{k}(u) -s = \frac{d}{dw} \hat F(G\cdot u;s)\, \ell_{k}(u) \cdot L\otimes G\, r_{k}(w)
\label{eq:almostlax1}
\end{equation}
Choosing $u=u_\pm$ and assuming $|u_+-u_-|$ small enough, 
the main assumption (\ref{eq:mainassump}) implies
 \begin{equation}\label{eq:convexred2NEW}
   \frac{d}{dw} \hat F(G\cdot u;s) \Bigr|_{u = u_\pm} = o(1)\qquad
    \textrm{as}\quad |u_-- u_+|\to 0.
\end{equation}
Differentiating (\ref{eq:almostlax1}) with respect to $u$ in the direction of $r_k(u)$, we obtain
\begin{align*}
 \nabla \lambda_k(u) \cdot r_k(u) = &\frac{d}{dw} \hat F(G \cdot u;s)  
 \nabla (\ell_{k}(u) \cdot L\otimes
              G r_{k}(u)) \cdot r_k(u) \\
              &+ \frac{d^2}{dw^2} \hat F(G \cdot u;s) 
               \bigl(\ell_{k}(u) \cdot L\bigr)\,\bigl(G\cdot r_{k}(u)\bigr)^2.
\end{align*}
Evaluating this relation at $u=u_\pm$ and taking in account (\ref{eq:convexred2NEW}),
we obtain
\begin{equation*}
 \frac{d^2}{dw^2} \hat F(G \cdot u_\pm ;s)
 =\frac{\nabla \lambda_k(u_\pm) \,\cdot\, r_k(u_\pm)}
 {\bigl(\ell_{k}(u) \cdot L\bigr)\,\bigl(G\cdot r_{k}(u)\bigr)^2}+o(1)
 \qquad \textrm{as}\quad |u_-- u_+|\to 1,
\end{equation*}
The conclusion follows from GNL condition (\ref{eq:GNL}) and assumption (\ref{eq:mainassump}).
\end{proof}

The previous results show that the existence of a $k$--radiative shock for 
(\ref{eq:system}) reduces to the study of a  scalar model of the form
   \begin{equation}
	   \begin{cases}
	       w_{t}+ \hat f(w)_{x} +q_{x}=0 & \\
	       -q_{xx} + R\,q +w_{x} =0, &
	   \end{cases}
	   \label{eq:scalar2}
       \end{equation}
provided (\ref{eq:mainassump}) is satisfied for the original shock $(u_{-},u_{+};s)$. 
In addition, if (\ref{eq:GNL})
    is also verified, we can assume, without loss of generality, 
    that the flux in (\ref{eq:scalar2}) is strictly convex 
    (see Proposition \ref{prop:convexreduction}). 
    In this framework, the reduction of (\ref{eq:system}) to  (\ref{eq:scalar2}) is 
    given by (\ref{eq:reduction}) and Lemma \ref{lem:reduction}, namely
	\begin{align*}
	     & w = G\cdot u, &  & \hat F(G\cdot u;s) = Q ( f(u) -  
	     f(u_{\pm}) - s(u-u_{\pm})  ), 
	\end{align*}
   with  $Q$ (and $P$ below) as before and, 
   adding if necessary a linear function to the flux,
   \begin{equation*}
       \hat F(w;s) = \hat f(w) - \hat f(w_{\pm}) - s (w-w_{\pm}).
   \end{equation*}
For the sake of clarity, we start by considering the GNL case and 
we postpone the general case at the end of the section.

If the flux in (\ref{eq:scalar2}) is strictly convex, 
Theorem \ref{theo:existence} guarantees the existence of a radiative shock 
   for that model, with at most a jump discontinuity, 
   which is indeed an admissible shock for the inviscid related 
   conservation law. To conclude with the results stated in Theorem 
   \ref{theo:mainintro}, we
   shall analyze that  discontinuity in  the corresponding
   radiative shock of the original vectorial case (\ref{eq:system}), 
   showing it forms 
   an admissible radiative shock for that system.

\begin{proof}[Proof of Theorem \ref{theo:mainintro}]
Let $(u_{-},u_{+};s)$ be an admissible $k$--shock for (\ref{eq:system}) and,  let us consider 
$w_{_{\pm}} = G\cdot u_{\pm}$.
Let $\mathcal{U}$, $\mathcal{W} = \Phi(\mathcal{U};s)$ be the neighborhoods of $u_{\pm}$ 
and $w_{\pm}$ given by Lemma \ref{lem:reduction} and Proposition \ref{prop:convexreduction}, 
and such that
	  \begin{equation}
	      \ell_{k}(u)\cdot L\otimes G\,r_{k}(u) >0
	      \label{eq:mainenlarged}
	  \end{equation}
	  for any $u\in\mathcal{U}$. 
    We start by proving that $(w_{-},w_{+};s)$ is an admissible shock for
  the reduced scalar conservation law (with strictly convex flux)
    \begin{equation}
	w_{t}+\hat f(w)_{x}=0.
	\label{eq:reducedagain}
    \end{equation}
   Indeed, from 
    \begin{equation*}
	  f(u_{+}) -  f(u_{-}) - s(u_{+}-u_{-})=0,
    \end{equation*}
    we obtain in particular
    \begin{equation*}
	\hat F(w_{+};s) = Q( f(u_{+}) -  f(u_{-}) - s(u_{+}-u_{-})) =0,
    \end{equation*}
that is, the Rankine--Hugoniot conditions of (\ref{eq:reducedagain}) 
for the shock $(w_{-},w_{+};s)$. 

Proceeding as in the proof of Proposition \ref{prop:convexreduction}, we obtain 
(\ref{eq:almostlax1}), that is
\begin{equation}
  \lambda_{k}(u) -s = (\hat f'(G\cdot u) -s) \ell_{k}(u) \cdot L\otimes G\, r_{k}(w).
\label{eq:almostlax1bis}
\end{equation}
Therefore the sign of $\hat f'(G\cdot u) -s$ is given by the sign of $\lambda_{k}(u) -s$ for any 
$u\in\mathcal{U}$, in view of (\ref{eq:mainenlarged}).
Thus, using (\ref{eq:almostlax1bis}) for $u=u_{\pm}$ and taking into account the admissibility 
condition of $u_{\pm}$, $\lambda_{k}(u_{+})<s<\lambda_{k}(u_{-})$, we obtain 
\begin{equation*}
 \hat f'(w_{-})>s>\hat f'(w_{+}),
\end{equation*}
that is, the discontinuity $(w_{-},w_{+};s)$ is admissible for (\ref{eq:reducedagain}).
	  
At this point, let $w$ be the (unique up to space shifts) radiative profile of the reduced
model (\ref{eq:scalar2}) given by Theorem  \ref{theo:existence}. 
Then, the above results guarantee the existence of the (unique up to space shift)  
radiative profile $(u,q)$ for (\ref{eq:system}), if  $|u_{-} - u_{+}|$ is sufficiently small.
Therefore we only have to prove that, if $w$ is discontinuous, the corresponding 
discontinuity in $u$ defines an admissible $k$--shock wave for the inviscid hyperbolic 
system of conservation laws $u_{t}+f(u)_{x}=0$. 
We shall perform this task as before, connecting the properties of the  shock for the 
reduced system with the ones of the shock of the original systems.
    
Denoting with $(w_{l},w_{r};s)$ the discontinuity of the radiating shock for the scalar reduced 
model, we have $w_{l,r}\in \mathcal{W}$ and therefore there exist unique $u_{l,r}\in\mathcal{U}$ 
solutions of
       \begin{equation}
	   \begin{cases}
	       G\cdot u = w_{l,r} &   \\
	       P(f(u) - f(u_{\pm}) - s(u-u_{\pm})  ) =0, & 
	   \end{cases}
	   \label{eq:systemred}
       \end{equation}
which are defined through the function $\Phi$ constructed in Lemma \ref{lem:reduction}. 
       
As before, let us start by showing the Rankine--Hugoniot condition for $(u_{l},u_{r};s)$,  namely
       \begin{equation*}
	     f(u_{l}) -   f(u_{r}) = s(u_{l} - u_{r}).
       \end{equation*}
Clearly, the above relation is equivalent to
       \begin{equation}
	     F(u_{l};s) =   F(u_{r};s).
	   \label{eq:RHrewritten}
       \end{equation}
 In order to prove (\ref{eq:RHrewritten}), we only have to prove 
       \begin{equation}
	   Q  F(u_{l};s) = Q  F(u_{r};s),
	   \label{eq:withQ}
       \end{equation}
because $P  F(u_{l};s) =0= P   F(u_{r};s)$ is given by  (\ref{eq:systemred})$_{2}$. 
Moreover, $Q  F(u;s) = \hat F(G\cdot u=w;s)$ for any $u\in\mathcal{U}$, namely, for 
$u=\Phi(w;s)$ and therefore (\ref{eq:withQ}) is precisely the Rankine--Hugoniot 
condition 
	   \begin{equation*}
	       \hat f(w_{l}) - \hat f(w_{r}) =s(w_{l}-w_{r})
	   \end{equation*}
for the reduced scalar model (\ref{eq:scalar2}).
       
We turn now to the proof of the Lax conditions for the $k$--shock $(u_{l},u_{r};s)$ 
of $u_{t} +   f(u)_{x}=0$, namely
       \begin{align}
	    &  \lambda_{k}(u_{r}) < s < \lambda_{k}(u_{l}), 
	   \label{eq:lax1}  \\
	    & \lambda_{k-1}(u_{l}) < s < \lambda_{k+1}(u_{r}).
	   \label{eq:lax2}
       \end{align} 
       Since we are dealing with weak shocks and the system is assumed 
       to be strictly hyperbolic, (\ref{eq:lax2}) follows from (\ref{eq:lax1}).
       Finally, using once again (\ref{eq:almostlax1}), this time  for 
       $u=u_{l,r}$, and taking into account the admissibility condition
       of $w_{l,r}= G\cdot u_{l,r}$, that is
       \begin{equation*}
	   \hat f'(w_{l})>s>\hat f'(w_{r}),
       \end{equation*}
       we obtain (\ref{eq:lax1}).
\end{proof}

\begin{proof}[Proof of Theorem \ref{theo:reg}]
    The proof is an immediate consequence
     of Proposition \ref{prop:reg} and Theorem \ref{theo:mainintro}.
\end{proof}

Let us pass now to the general case, namely when condition (\ref{eq:GNL}) is violated. 
As we have shown in Section \ref{subsec:nonconvex}, the existence of radiative shock without 
convexity assumptions is guaranteed only when the radiative effect is sufficiently dissipative. 
Hence, we shall consider the following system
\begin{equation}
    \begin{cases}
	u_{t}+ f(u)_{x} +Lq_{x}=0 & \\
	-\e\, q_{xx} + R\,q +G\cdot u_{x} =0, &
    \end{cases}
    \label{eq:system2}
\end{equation}
with $0<\e \ll 1$. It is worth to observe that the smallness in $\e$ 
is needed only for the existence of the profile for the reduced scalar  model
    \begin{equation}
	       \begin{cases}
		   w_{t}+ \hat f(w)_{x} +q_{x}=0 & \\
		   -\e\, q_{xx} + R\,q +w_{x} =0 &
	       \end{cases}
	       \label{eq:scalar3}
	   \end{equation}
and it does not play any role in the connections between the 
admissibility conditions for the jumps of that model and the ones of (\ref{eq:system2}).

\begin{theorem}\label{theo:RHgeneral}
Let $(u_{-},u_{+};s)$ be an admissible $k$--shock for (\ref{eq:reduced}) and
assume (\ref{eq:mainassump})  holds. 
 If $|u_{-} - u_{+}|$  and $\e$ are sufficiently small, then there exists a (unique up to shift) 
 admissible radiative shock wave $(u=u(x-st),q=q(x-st))$ of (\ref{eq:system2}) such that 
 $(u(\pm\infty),q(\pm\infty)) = (u_{\pm},0)$.
\end{theorem}

\begin{proof}
    The proof of this theorem follows the same lines of the one of 
 Theorem \ref{theo:mainintro}. In particular, the existence of a 
    radiative shock for (\ref{eq:system2}) comes from Theorem 
    \ref{theo:generalflux} and the reduction   to 
    (\ref{eq:scalar3}). Clearly, the analysis of Rankine--Hugonoit 
    conditions is made as before, because it is independent from 
    convexity assumptions.
    Theorefore we are left to the proof that, given the neighborhoods 
    $\mathcal{U}$ and $\mathcal{W}$ as before,  a shock 
    $(u_{l},u_{r};s)$, $u_{l,r}\in\mathcal{U}$, verifies the Liu 
    E--condition for (\ref{eq:reduced}) if and only if the 
    corresponding shock $(w_{l},w_{r};s)$, $w_{l,r}\in\mathcal{W}$, 
    verifies the Oleinik condition for 
    \begin{equation}
	w_{t}+ \hat f(w)_{x}=0.
        \label{eq:reducedlast}
    \end{equation}
 Let  $(u_{-}, u_{+};s)$ be a $k$-shock of (\ref{eq:system}) which 
    verifies the
    (strict) Liu E--condition: denoting with $u_{k}(\tau)$ the 
    $k$--th shock curve, $u_{-} = u_{k}(0)$, $u_{+} = u_{k}(\bar 
    \tau)$, then 
    \begin{equation}
        s = s_{k}(\bar \tau) < s_{k}(\tau)
        \label{eq:liu}
    \end{equation}
    for any $\tau$ between $0$ and $\bar \tau$. As before, consider 
    $w_{\pm} = G\cdot u_{\pm}$ and assume $w_{-}>w_{+}$. By 
    definition of shock curve, we have
    \begin{align*}
        0 &= f(u_{k}(\tau)) - f(u_{-}) - s_{k}(\tau)(u_{k}(\tau) - u_{-}) \\
	&= F(u_{k}(\tau);s) + (s - s_{k}(\tau))(u_{k}(\tau) - u_{-}) \\
	&= \hat F(G\cdot u_{k}(\tau);s) L + (s - s_{k}(\tau))(u_{k}(\tau) - u_{-})
    \end{align*}
    which implies
    \begin{equation}
	\hat F(G\cdot u_{k}(\tau);s) L = (s_{k}(\tau) - s)(u_{k}(\tau) - u_{-}).
        \label{eq:liu2}
    \end{equation}
    Since for small shocks, that is $|\tau|$ small, $u_{k}(\tau) - 
    u_{-} = \tau r_{k}(u_{-}) + o(\tau)$, form (\ref{eq:liu2}) we 
    conclude
    \begin{equation*}
	\hat F(G\cdot u_{k}(\tau);s) L = (\tau\,r_{k}(u_{-}) + o(\tau))(s_{k}(\tau) - s)
    \end{equation*}
    and, multiplying that relation on the left for $\ell_{k}(u_{-})$
    \begin{equation}
	\hat F(G\cdot u_{k}(\tau);s)\ell_{k}(u_{-})\cdot L = (\tau + o(\tau))(s_{k}(\tau) - s).
        \label{eq:liu3}
    \end{equation}
    Moreover, we differentiate $w_{k}(\tau) := G\cdot u_{k}(\tau)$ 
    with respect to $\tau$ to conclude $\dot{w}_{k}(\tau) = G\cdot 
    \dot{u}_{k}(\tau) = G\cdot r_{k}(u_{-}) + o(\tau)$, and, since we 
    are dealing with small shocks and $G\cdot r_{k}(u_{-}) \neq 0$ 
    for (\ref{eq:mainassump}), 
    it follows that $w_{k}(\tau)$ is  decreasing for $\tau$ between 
    $0$ and $\bar\tau$. Hence, if $G\cdot r_{k}(u_{-}) > 0$ (resp. 
    $G\cdot r_{k}(u_{-}) < 0$), then $\bar\tau <0$ (resp.\ $\bar\tau 
    >0$) and $\ell_{k}(u_{-})\cdot L >0$ (resp. $\ell_{k}(u_{-})\cdot 
    L <0$) using again (\ref{eq:mainassump}). Therefore, for $|\tau|$ 
    small and between $0$ and $\bar\tau$, $\tau + 
    o(\tau) < 0$ (resp.\ $\tau + 
    o(\tau) >0$) if  $\ell_{k}(u_{-})\cdot L >0$ (resp. $\ell_{k}(u_{-})\cdot 
    L <0$) and therefore from (\ref{eq:liu}) and (\ref{eq:liu3}) we conclude
    $\hat F(G\cdot u_{k}(\tau);s)<0$ for any $\tau$ between $0$ and 
    $\bar\tau$, that is, the strict Oleinik 
    condition for the shock $(w_{-},w_{+};s)$ of 
    (\ref{eq:reducedlast}).
    
    Let us now consider a radiative profile for (\ref{eq:scalar3}) 
    which has a discontinuity $(w_{l},w_{r};s)$ admissible 
    for (\ref{eq:reducedlast}). As in the proof 
    of Theorem \ref{theo:mainintro} for the GNL case, there exist 
    unique $u_{l,r}$ in the neighborhood under consideration such 
    that the discontinuity $(u_{l},u_{r};s)$ verifies the 
    Rankine--Hugoniot conditions 
    of (\ref{eq:reduced}) and in 
    particular it belongs to a $k$--shock curve $u_{k}(\tau)$ with 
    $u_{k}(0)=u_{l}$, $u_{k}(\tau_{r}) = u_{r}$.
    The relation 
    \begin{equation*}
	f(u_{k}(\tau)) - f(u_{l}) - s_{k}(\tau)(u_{k}(\tau) - 
		u_{l}) = 0
    \end{equation*}
    implies this time
    \begin{equation*}
	(\hat F(G\cdot u_{k}(\tau);s)- 
	\hat F(G\cdot u_{l,r}(\tau);s) L = (s_{k}(\tau) - s)(u_{k}(\tau) - 
		u_{-}),
    \end{equation*}
    which becomes for $|\tau|$ small
    \begin{equation*}
	(\hat F(G\cdot u_{k}(\tau);s) - 
	\hat F(G\cdot u_{l,r}(\tau);s) )\ell_{k}(u_{l})\cdot L = 
	(s_{k}(\tau) - s)(\tau + o(\tau) ).
    \end{equation*}
   Moreover, using (\ref{eq:mainenlarged}) this time for $u_{l}$, we 
   argue as before to obtain the 
   (strict) Liu E--condition
   \begin{equation*}
       s=s_{k}(\tau_{r}) < s_{k}(\tau)
   \end{equation*}
   for any $\tau$ between $0$ and $\tau_{r}$ from the (strict) 
   Oleinik condition
   \begin{equation*}
       \hat F(G\cdot u_{k}(\tau);s) < 
	       \hat F(G\cdot u_{l,r}(\tau);s)
   \end{equation*}
   for any $w_{k}(\tau) = G\cdot u_{k}(\tau)\in (w_{r},w_{l})$
   and the proof is complete.
\end{proof}

\section*{Acknowledgements}
The research of the authors was partially supported by the European
IHP project 'HYKE', contract \# HPRN-CT-2002-00282. 
It was achieved in part when the third author was a Visiting Professor 
INdAM-GNAMPA at the Dipartimento di Matematica Pura e 
di Matematica Pura e Applicata, Universit`a degli Studi dell'Aquila.

\providecommand{\bysame}{\leavevmode\hbox to3em{\hrulefill}\thinspace}
\providecommand{\MR}{\relax\ifhmode\unskip\space\fi MR }
\providecommand{\MRhref}[2]{%
  \href{http://www.ams.org/mathscinet-getitem?mr=#1}{#2}
}
\providecommand{\href}[2]{#2}

\end{document}